\documentclass[]{svjour3}       
\smartqed 
\RequirePackage{fix-cm}
\usepackage[utf8]{inputenc}
 \usepackage{float}
 \usepackage{varioref}                             
 \usepackage{wrapfig}                              
 \usepackage{threeparttable}                       
 \usepackage{dcolumn}   
   \newcolumntype{d}{D{.}{.}{-1}}
 \usepackage{graphicx}
 \usepackage{epstopdf}
 \usepackage{subfigure}                            
 \usepackage{subfigmat}                            
 \usepackage{fancyvrb}                             
    \fvset{fontsize=\footnotesize,xleftmargin=2em}
\usepackage{lettrine}                             
\usepackage[colorlinks]{hyperref}                 
\usepackage{amsmath}
\usepackage{longtable,tabularx}
\usepackage{placeins}
\usepackage{multirow}
\usepackage{makeidx}
\usepackage{mdframed}
\usepackage{xcolor}
\usepackage{algorithm}
\usepackage[noend]{algpseudocode}
\makeatletter
\def\BState{\State\hskip-\ALG@thistlm}
\makeatother
\newlength\tindent
\newcommand{\pd}[2] {\frac{\partial #1}{\partial #2}}
\newcommand{\td}[2] {\frac{d #1}{d #2}}

\newcommand{\norm}[1]{\left\lVert#1\right\rVert}

\journalname{Numerical Algorithms}
\begin{document}
\title{Approximate Linearization of Fixed Point Iterations}
\subtitle{Error Analysis of Tangent and Adjoint Problems Linearized about Non-Stationary Points}
\author{Emmett Padway \and Dimitri Mavriplis}
\institute{E. Padway \at
              National Institute of Aerospace, 100 Exploration Way, Hampton, VA, 23666 \\
              \email{emmett.padway@nianet.org}           
           \and
           D. Mavriplis \at
              University of Wyoming, Department of Mechanical Engineering, 1000 E. University Ave, Laramie, WY, 82071}

\date{Received: date / Accepted: date}



\maketitle

\begin{abstract}
Previous papers have shown the impact of partial convergence of discretized PDEs
on the accuracy of tangent and adjoint linearizations. A series of papers
suggested linearization of the fixed point iteration used in the solution
process as a means of computing the sensitivities rather than linearizing the
discretized PDE, as the lack of convergence of the nonlinear problem indicates
that the discretized form of the governing equations has not been satisfied.
These works showed that the accuracy of an approximate linearization depends in
part on the convergence of the nonlinear system. This work shows an error
analysis of the impact of the approximate linearization and the convergence of
the nonlinear problem for both the tangent and adjoint modes and provides a series
of results for an exact Newton solver, an inexact Newton solver, and a low
storage explicit Runge-Kutta scheme to confirm the error analyses.
\end{abstract}
\section{Introduction}
\lettrine[nindent=0pt]{T}{angent} and Adjoint methods have become a larger part
of the field of Computational Fluid Dynamics (CFD)  as computers and algorithms
have developed and a concurrent push for automatic design has strengthened. Both
of these methods are derived from the requirement that the discretized form of
the governing Partial Differential Equations (PDE) is satisfied; for steady state problems this is characterized
by the magnitude of the residual vector of the spatially discretized terms approaching machine zero. The adjoint linearization (either through the
discrete or continuous methods) is more commonly used as it scales independent
of the number of design variables and is used for error estimation
\cite{MAMN,DDAMR}. These methods are highly accurate and are often verified
through use of the complex-step finite-difference method as this does not suffer
from round off error \cite{SKN_Thesis}. However, for cases where the discretized
governing equations are not satisfied, the complex-step finite-difference method gives the
sensitivity of the solution process and both the adjoint and tangent linearizations are no longer
exact linearizations with the adjoint losing duality to the primal problem. Such
cases have been shown to be very common in CFD simulations
\cite{JKDD,JKQW,PTA-devel} and in such cases the sensitivities suffer greatly in
terms of accuracy and exhibit a high dependence on the state of the simulation
at termination. 

To address these cases Padway and Mavriplis \cite
{PTA-devel,PTA-iqn} investigated linearizing the fixed point iteration used to
solve the nonlinear problem and obtained the sensitivities of the solution
process. Through this method they showed better sensitivity behavior and
successful optimizations; they also applied this method to error estimation and
adaptive mesh refinement
\cite{PTA-amr}. This method inherits some of the characteristics shown in the
 work on one-shot optimization methods \cite{OneShotOpt} including the
 guaranteed convergence of the tangent and adjoint problems, due to their
 similarity to the black box approach \cite{BB_Sagebaum}. Padway and
 Mavriplis \cite{PTA-iqn} showed by numerical experiment that for an
 approximately linearized quasi-Newton fixed point iteration, convergence of
 the non-linear problem led to a decrease in error from the approximate
 linearization. The issue of inexact linearization of a linear system solve has
 been investigated previously to produce consistent automatic differentiation
 of linear system solves in segregated solvers \cite{JDM_LinAlg}. These
 linearizations would be necessary in the "piggy-back" iterations of the
 one-shot adjoint method \cite{OSA_Gunther} when applied to implicit nonlinear
 solvers. In this paper the authors show through a series of proofs and
 numerical experiments, that -- providing the residual operator is linearized
 exactly -- the error from inexact linearization of the fixed point iteration
 approaches 0 at the same rate of the nonlinear problem convergence, a necessary condition for approximate linearizations to be useful in
 applications. Additionally,  the authors show that for
 nonconverged problems, the error in the sensitivity has two terms, one of
 which scales with the residual, and the other of which scales with the error
 in the approximate linearization; the second of which is multiplied by the
 derivative of the contractive fixed point iteration; a
 discovery that allows a user to select the level of error they desire in their
 sensitivities. Additionally, this allows users to know that by converging the
 nonlinear problem they can obtain greater accuracy in the sensitivities
 despite the inaccuracy of the approximate linearization of the fixed point
 iteration. As problems with nonconvergent fixed points have become more
 prevalent as the field has attempted more difficult problems, methods like
 those mentioned above show more promise. In the realm of design optimization,
 which is the focus of this paper, inaccurate adjoints can lead to inaccurate
 sensitivities, which can change the course of the design cycle and lead to
 stagnation; as the Karush-Kuhn-Tucker (KKT) conditions, which govern
 termination, require that the gradient vanish at a local extremum \cite
 {SNopt2}. As such having a good estimate on the error of the sensitivities for
 such a method (when compared to the complex-step finite-difference method)
 that can provide better designs than those provided by the solution of the
 classical steady-state adjoint is highly desirable. In a similar vein, Brown
 and Nadarajah \cite{SKN} investigate a bound for the error in the adjoint
 computed sensitivities arising from partial convergence of the primal problem
 for problems that show smooth convergence behavior.  

This paper is organized as follows. We first consider the details of the CFD
solver used in this work, explaining the discretization, the solution
methodologies, and the classical steady-state adjoint and tangent problems. Then
we show the pseudo-time accurate tangent and adjoint linearizations to obtain
design sensitivities along with a discussion of the approximation decisions.
After which we show analyses of the error in the design sensitivities of both
linearizations. We then show numerical results for sensitivity computations for
an airfoil in transonic flow, thus confirming the analyses and showing that
these analyses hold for both explicit and implicit solvers. 

\section{Background and In-House Solver}

\subsection{Governing Equations}
We developed an in-house flow solver to solve the steady-state Euler equations
on unstructured meshes. The steady-state compressible Euler equations (which may
also be referred to as the primal or analysis problem) can be written as follows:
\begin{equation} \label{eq:DivFlux}
    \nabla \cdot F(u(D)) = 0,
\end{equation}
where $u$ is the conservative variable vector, $D$ is the design variable vector, and $F(u)$ is the conservative variable flux. Equation (\ref{eq:DivFlux}) can also be written as:
\begin{equation} \label{eq:Res}
R(u(D), D) = 0,
\end{equation}
where R is the residual operator. 
\subsection{Spatial Discretization}
The residual about the closed control volume is given as
\begin{equation}
    R = \int \left[F(u(D)) \right] \cdot n(x(D)) \mathrm{dB} = \sum_{i=1}^{n_{edge}}
F_{e_i}^{\bot}(u(D), n_{e_i}(x(D)))B_{e_i}(x(D)),
\end{equation} 
 this is the operator in the aforementioned requirement for the
 adjoint and tangent systems. In the discretized form of the
 residual operator, $e_i$ is a given edge of the triangulation, $n_{edge}$ is
 the number of edges in the triangulation, $B$ is the element boundary, $x$ is
 the mesh point coordinate vector, $F$ is the numerical flux across the element
 boundary, and $n$ is the edge normal on the element boundary. The solver used
 in this work is a steady-state finite-volume cell-centered Euler solver with
 second-order spatial accuracy for triangular elements. Second-order accuracy
 is implemented through weighted least squares gradient reconstruction
\cite{LSQR}. The solver has three different flux calculations implemented and
 linearized, these are the Lax-Friedrichs\cite{LF}, Roe\cite{Roe}, and Van
 Leer\cite{VL} schemes. In this work, only the Van-Leer scheme
 is used due to the suitability of the results at a decreased expense and
 increased stability when compared to the Roe flux.

\FloatBarrier
\subsection{Steady-State Solver}
To solve the discretized PDE shown above we define a fixed point iteration (N): 
\begin{equation}
\begin{aligned}
u^{k} &= N(u^{k-1}, x(D)) = u^{k-1} + H(u^{k-1}, x(D)) \\
&= u^{k-1} + A(u^{k-1},x(D))R(u^{k-1},x(D)),	
\end{aligned}
\end{equation}
where $H$ is the solution increment (also referred to as $\Delta u$), $A$ is an operator that defines the nonlinear solver and $x$ is the vector of mesh point coordinates in the triangulation which is a function of the design variable vector (often a parameterization) held fixed at each design iteration.
The solver technology for this code uses either explicit pseudo-time stepping
with a forward Euler scheme or a low storage five stage Runge-Kutta scheme, or
implicit pseudo-time stepping with a quasi-Newton method. The quasi-Newton
method is implemented using pseudo-transient continuation (PTC) with a BDF1
pseudo temporal discretization scheme. For Newton's method the time-stepping
procedure is written as:

\begin{equation} \label{eq:Nts}
u^{k} = u^{k-1} + \Delta u,
\end{equation}
where we compute $\Delta u$ by solving the following system of linear equations.
\begin{equation} \label{eq:NM}
\left[P \right] \Delta u = -R(u),
\end{equation}
where $\left[P\right]$ denotes a left hand matrix that determines the nonlinear solver. If $P$ is the exact linearization of the residual operator this is an exact Newton solver; otherwise we recover an inexact-quasi-Newton solver of which the most typical matrix is the linearization of a lower order discretization with a pseudo-time term. 
We can substitute the expression for $\Delta u$ into the time-stepping equation
(\ref{eq:Nts}) to obtain the final form of this equation
\begin{equation} \label{eq:PTCBDF1}
\begin{aligned}
u^{k} &= u^{k-1} - \left[P_{k-1}\right]^{-1}R.
\end{aligned}
\end{equation} 
Here $\left[P_{k-1}\right]$ is the Jacobian of the first order spatial discretization
augmented with a diagonal term to ensure that it is diagonally dominant, shown
in equation (\ref{eq:PC1}).  The use of the Jacobian of the first-order spatial discretization turns this into an inexact-quasi-Newton solver, which is often used due to the robustness and ease of implementation when compared to the Jacobian of a second-order discretization. This is especially desirable in this work because the inspiration behind the methods that are analyzed herein are methods used for optimization or error esimation where the flow solver does not converge and the added robustness is desired.
\begin{equation} \label{eq:PC1}
\left[P_{k-1}\right] = \left[\pd{R}{u^{k-1}}\right]_1 + \frac{vol}{\Delta t CFL} ,
\end{equation}
where $vol$ is the volume (or area) of the element, $CFL$ is the CFL number at the given iteration and $\Delta t$ is a local time step given as:
\begin{equation}
    \Delta t_i = \frac{r_i}{\sqrt{(u^2 + v^2)} + c} . 
\end{equation}
Here $r_i$ is the circumference of the inscribed circle for mesh cell $i$, $u$ and
$v$ are the horizontal and vertical velocity components respectively, and $c$ is the
speed of sound in the triangular element. Please note the subscript on the
Jacobian in equation (\ref{eq:PC1}) denotes that it is the Jacobian of the first
order accurate spatial discretization; the subscript of 2 would indicate the
Jacobian of the second order accurate discretization, and that for the above
definition of the $A$ operator this would imply:
\begin{equation}
A(u^{k-1},D) = -\left[P_{k-1}\right]^{-1}.
\end{equation}

Furthermore, the CFL is scaled either with a simple ramping coefficient ($\beta$) and cap criterion:
\begin{equation}
   CFL^{k+1} = min(\beta \cdot CFL, CFL_{max})
\end{equation}
or with a line search and CFL controller \cite{DJMLS,BALS}, which seeks to minimize the $L_2$ norm of the pseudo-temporal residual, defined as:
\begin{equation}
R_t(u + \alpha\Delta u) = \frac{vol}{\Delta t}\alpha\Delta u^k + R(u + \alpha\Delta u^k).
\end{equation}
\FloatBarrier

\noindent In order to solve the linear system that arises in the quasi-Newton
method we use a flexible GMRES solver \cite{GMRES} preconditioned by 
Gauss-Seidel iterations or the Gauss-Seidel iterations alone. The flexible GMRES solver is also used to solve the stiff steady state tangent
and adjoint systems as well.

\subsection{A Review of Tangent and Adjoint Systems}
\subsubsection{Tangent Formulation}
For an aerodynamic optimization problem, we consider an objective functional
$L(u(D),x(D))$, for example lift or drag. In order to obtain an expression
for the sensitivities we take the derivative of the objective functional
\cite{mavriplis_vki}:
\begin{equation} \label{eq:tan}
\frac{dL}{dD} = \pd{L}{x}\td{x}{D}+\pd{L}{u}\td{u}{D}.
\end{equation} 
For the above expression $\pd{L}{x}$ and $\pd{L}{u}$ can be directly obtained
by differentiating the corresponding subroutines in the code and $\td{x}{D}$ is
dependent on the choice of mesh motion scheme and design variables. It is not
possible to obtain $\td{u}{D}$ through linearization of the subroutines in the
code without linearizing the entire analysis solution process, as will be
covered in later sections. In order to solve for this term we use the
constraint that for a fully converged flow $R(u(D), x(D)) = 0$. By taking the
derivative of the residual operator we obtain the equation below:

\begin{equation}
\left[\pd{R}{x}\right]\td{x}{D} + \left[\pd{R}{u}\right]_2\td{u}{D} = 0 .
\end{equation}

We can isolate the sensitivity of the residual to the design variables to obtain the tangent system:
\begin{equation} \label{eq:tanlin}
\left[\pd{R}{u}\right]_2\td{u}{D} = -\left[\pd{R}{x}\right]\td{x}{D}.
\end{equation}
We solve this linear system, using hand differentiated subroutines to provide
the left hand matrix $\left[\pd{R}{u}\right]_2$, the right hand side
$\left[\pd{R}{x}\right]\td{x}{D}$ (which scales with the design variables), and
obtain $\td{u}{D}$. We then substitute $\td{u}{D}$ into equation (\ref{eq:tan})
to obtain the final sensitivities.

\subsubsection{Discrete Adjoint Formulation}
The adjoint formulation begins with the same sensitivity equation:
\begin{equation} \label{eq:tanlinadj}
\frac{dL}{dD} = \pd{L}{x}\td{x}{D}+\pd{L}{u}\td{u}{D}.
\end{equation} 
Using the condition $R(u(D),x(D))=0$, we return to equation (\ref{eq:tanlin}) and
pre-multiply both sides of the equation by the inverse Jacobian matrix to
obtain:
\begin{equation}
\td{u}{D} = -\left[\pd{R}{u}\right]_2^{-1}\left[\pd{R}{x}\right]\td{x}{D}.
\end{equation}
Substituting the above expression into the sensitivity equation yields:
\begin{equation}
\frac{dL}{dD} = \pd{L}{x}\td{x}{D}-\pd{L}{u}\left[\pd{R}{u}\right]_2^{-1}\left[\pd{R}{x}\right]\td{x}{D}.
\end{equation} 
We then define an adjoint variable ${\bf \Lambda}$ such that:
\begin{equation}
{\bf \Lambda}^T = -\pd{L}{u}\left[\pd{R}{u}\right]_2^{-1},
\end{equation}
which gives an equation for the adjoint variable:
\begin{equation}
\left[\pd{R}{u}\right]_2^T{\bf \Lambda} = -\left[\pd{L}{u}\right]^T.
\end{equation}
We can solve this linear system and obtain the sensitivities for the objective function as follows:
\begin{equation}
\frac{dL}{dD} = \left[\pd{L}{x}+{\bf \Lambda}^T\pd{R}{x}\right]\td{x}{D}.
\end{equation} 
The adjoint system is of interest, because as mentioned in the introduction, it
results in an equation for the sensitivity that does not scale with the number
of design variables. 

\FloatBarrier
\section{Development of the Pseudo-Time Accurate Tangent} 
In this section we show the previously derived pseudo-time accurate tangent formulations
\cite{PTA-devel,PTA-iqn}. We also discuss the impacts on implementation and
accuracy of some of the assumptions made in the previous formulation; the assumptions are exact when the linear
system at each iteration is solved to machine precision. In such cases the
tangent sensitivities would exactly correspond to the sensitivities of solution
process itself, and would correspond exactly to the sensitivities obtained by the
complex-step finite difference method at each step of the solution process.  We note that in this section and
the derivation of the adjoint section, when taking the derivative of an operator
with respect to the design variables, $\pd{}{D}$ is an abbreviation of
$\pd{}{x}\td{x}{D}$.
\subsection{Tangent System for quasi-Newton Method}

Previous works on the pseudo-time accurate adjoint \cite{PTA-devel,PTA-iqn} have contained investigations on the desireability of an objective function averaged in pseudo-time over the last $m+1$ steps, expressed as:
\begin{equation}
L = L(u^n,u^{n-1},...,u^{n-m},D).
\end{equation}
The sensitivity equation for such a problem is
\begin{equation} \label{eq:TanSens}
\td{L}{D} = \pd{L}{D} + \pd{L}{u^n}\td{u^n}{D} + \pd{L}{u^{n-1}}\td{u^{n-1}}{D} + ... + \pd{L}{u^{n-m}}\td{u^{n-m}}{D}.
\end{equation}
In this work we focus on objective functions calculated only at the final iteration/pseudo-time step; such a decision greatly simplifies the analyses and they can still be straightforwardly extended to multiple iteration objective functions.  We recover
\begin{equation}
\td{L}{D} = \pd{L}{D} + \pd{L}{u^n}\td{u^n}{D},
\end{equation}
which is identical to the sensitivity equation of the steady state tangent system (\ref{eq:tan}).
For this section, we begin from equation (\ref{eq:PTCBDF1})
\begin{equation} \label{eq:PTCBDF2}
\begin{aligned}
u^{k} &= u^{k-1} - \left[P_{k-1}\right]^{-1}R(u^{k-1}, x(D)),
\end{aligned}
\end{equation} 
where  $\left[P_{k-1}\right]$ is again a first-order accurate Jacobian augmented with
a diagonal term to ensure that it is diagonally dominant, as described in
equation (\ref{eq:PC1}):
\begin{equation} \label{eq:qNte}
\left[P_{k-1}\right] = \left[\pd{R}{u^{k-1}}\right]_1 + \frac{vol}{\Delta t CFL} .
\end{equation}
If we take the derivative of each side of equation (\ref{eq:PTCBDF2}) we obtain:
\begin{equation} \label{eq:Full_dqNM}
\begin{aligned}
\frac{du}{dD}^{k} &= \frac{du^{k-1}}{dD} - \left[P_{k-1}\right]^{-1}\left[ \pd{R}{D} + \left[\pd{R}{u^{k-1}}\right]_2\frac{du}{dD}^k\right] \\
&- \left[ \td{\left[P_{k-1}\right]^{-1}}{D}\right]R(u^{k-1}, x(D)) .
\end{aligned}
\end{equation}
Here, rather than neglecting the change in the preconditioner we use a
definition of the derivative of the matrix inverse defined by:
\begin{equation} \label{eq:dKinv}
\td{\left[K\right]^{-1}}{x} = -\left[K\right]^{-1}\left[ \td{K}{x}\right] \left[K\right]^{-1},
\end{equation}
 as applied in previous work \cite{PTA-iqn}. This
 assumption will not be exact for any case in which the linear system solution
 is not exact to machine precision. However, it will allow us to have a clearly
 defined source of error in a computation and evaluate how much the linear
 tolerance of the system affects the sensitivity computation. Such an investigation was not done in the previous work, but is the
 thrust of this paper, and the effect of this approximation will be examined,
 as well as an investigation into the behavior of other approximate
 linearizations. By computing the total derivative of the nonlinear solver--
 as in equation(\ref{eq:Full_dqNM})-- and  applying (\ref
 {eq:dKinv}) with $K = P_{k-1}$ into equation (\ref{eq:Full_dqNM}), the below
 expression is generated: 
\begin{equation}
\begin{aligned}
\td{u^k}{D} &= \td{u^{k-1}}{D} - \left[ P_{k-1}\right]^{-1}\left[\left[\pd{R}{u^{k-1}}\right]_2\td{u^{k-1}}{D} + \left[\pd{R}{x}\right]\td{x}{D}\right] \\
&+ \left[ P_{k-1}\right]^{-1}\td{\left[ P_{k-1}\right]}{D}\left[ P_{k-1}\right]^{-1}{R(u^{k-1},x(D))}.
\end{aligned}
\end{equation}
Expanding the total derivative shows that this is in fact the sum of two
matrix vector products:
\begin{equation} \label{eq:IQN2_TanR}
\begin{aligned}
\td{u^k}{D} &= \td{u^{k-1}}{D} - \left[
P_{k-1}\right]^{-1}\left[\left[\pd{R}{u^{k-1}}\right]_2\td{u^{k-1}}{D} + \left[\pd{R}{x}\right]\td{x}{D}\right]  \\
&+\left[ P_{k-1}\right]^{-1}\left[\pd{P_{k-1}}{u^{k-1}}\td{u^{k-1}}{D} + \pd{P_{k-1}}{x}\td{x}{D}\right] \left[ P_{k-1}\right]^{-1}{R(u^{k-1},x(D))} .
\end{aligned}
\end{equation}
While the full hessian is an undesireable item to compute, this formulation
requires two hessian vector products that can be obtained through complex
perturbations to the conservative variables and nodal coordinate vectors, and
subsequent evaluation of the Jacobian. This use of Frech\'et differentiation obviates the need for the full
hessian computation. Furthermore, one of the matrix inverses can be removed by
reusing the computation of $\Delta u$ in (\ref{eq:PTCBDF1}) and rewriting the equation as:
\begin{equation}
\begin{aligned}
\td{u^k}{D} &= \td{u^{k-1}}{D} - \left[
P_{k-1}\right]^{-1}\left[\left[\pd{R}{u^{k-1}}\right]_2\td{u^{k-1}}{D} + \left[\pd{R}{x}\right]\td{x}{D}\right]  \\
&+\left[ P_{k-1}\right]^{-1}\left[\pd{P_{k-1}}{u^{k-1}}\td{u^{k-1}}{D} + \pd{P_{k-1}}{x}\td{x}{D}\right] \Delta u .
\end{aligned}
\end{equation}

This can then be rewritten as:
\begin{equation}
\frac{du^k}{dD} = \frac{du^{k-1}}{dD} + \Delta\left(\td{u}{D}\right) ,
\end{equation}
where $\Delta \td{u}{D}$ is solved for in the following linear system:
\begin{equation} \label{eq:qNM_ddu}
\begin{aligned}
\left[P_{k-1}\right]\Delta\left(\td{u}{D}\right) &=
- \left[\left[\pd{R}{u^{k-1}}\right]_2\td{u^{k-1}}{D} + \left[\pd{R}{x}\right]\td{x}{D}\right] \\
  &+ \left[\pd{P_{k-1}}{u^{k-1}}\td{u^{k-1}}{D} + \pd{P_{k-1}}{x}\td{x}{D}\right]
\Delta u .
\end{aligned}
\end{equation} 
 It is important to note again that this expression is only exact
 for machine-zero solution of the linear system at every iteration, but for
 those cases we will have exact correspondence between this and the
 complex-step finite-difference computed sensitivities. This saves us from
 needing to differentiate the entire linear solution process, which is
 intractable for many cases for the forward mode, and even more difficult for
 the reverse mode. The reverse differentiation of a Krylov solver would be an
 onerous task that would yield little gain. One additional point is that with
 this method there are no conditions on exact duals of the linear solver, and
 one can use different solvers for the forward and the reverse linearizations.
 We can also note that, where in the initial formulation-- in equation (\ref
 {eq:IQN2_TanR})-- we see 3 approximate linear solves; when we group terms and
 use the already stored information from the nonlinear solution process we are
 left with only one approximate linear solve, the same as in the analysis
 problem's nonlinear solver. Although we are not limited to using the same
 linear solver for the analysis and tangent problems and its dual for the
 adjoint problem, by doing so we gain in that we are guaranteed to have a
 converging linear solver in both the tangent and adjoint modes. We would then
 be solving using the same algorithm that ran without divergence through the
 analysis portion, when we apply this algorithm to the tangent and adjoint
 problem we have the same eigenvalues, same condition numbers, and same
 convergence properties and the linear solvers will behave similarly.
 If the analysis problem linear solver did not diverge, neither will
 the linear solvers of the tangent and the adjoint problems.

\section{Development of the Pesudo-Time Accurate Adjoint}
The pseudo-time accurate adjoint method is drawn from the
derivation of the unsteady adjoint; we look at each pseudo-time
step and work backwards through pseudo-time to get the pseudo-time accurate
adjoint solution. In this derivation the objective function is a pseudo-time
averaged functional, averaged over the last  $m+1$ steps for a program that runs
through n pseudo-time steps.
From this we define an objective function L as:
\begin{equation}
L = L(u^n,u^{n-1},...,u^{n-m},D) ,
\end{equation} 
where $u^n$ is the conservative variable vector at the final time step n and D
is the design variable vector. For the constraint we cannot select $R(u,D) = 0$,
as this is not true at each pseudo-time step; instead, we select a constraint
based on the pseudo-time evolution of the solution, for which, the $k^{th}$ constraint
will be referred to as $G^{k}$. We know because we are using first order time-stepping,
that the constraint is dependent only on the old time-step, the new time-step, and
the design variables, expressed as follows:
\begin{equation}
G^k = G^k(u^k, u^{k-1}, D) = 0 .
\end{equation}
\\
We define an augmented objective function with n constraints and n Lagrange multipliers:
\begin{equation} \label{eq:AugObj}
\begin{aligned}
J(D, u^n, u^{n-1}, u^{n-2},..., \Lambda ^n, \Lambda^{n-1}, ...) &= L(u^n,u^{n-1}, ..., u^{n-m}, D) \\
&+ \Lambda ^{nT}G^n(u^n(D), u^{n-1}(D),D) \\
&+ ... \\ 
&+ \Lambda ^{1T}G^1(u^1(D), u^0(D),D) .
\end{aligned}
\end{equation}
\\
In order to get an expression for the adjoint we take the derivative of the
augmented objective function with respect to the conservative variables at
different pseudo-time steps, and choose Lagrange multipliers such that these
partial derivatives are equal to 0:
\begin{equation}
\begin{aligned}
\pd{J}{u^{n}}  &= \pd{L}{u^{n}} + \Lambda^{nT}\pd{G^{n}}{u^{n}} = 0 ,\\
\pd{J}{u^{n-1}}&= \pd{L}{u^{n-1}} + \Lambda^{nT}\pd{G^{n}}{u^{n-1}} + \Lambda^{n-1T}\pd{G^{n-1}}{u^{n-1}} = 0, \\
\pd{J}{u^{n-2}}&= \pd{L}{u^{n-2}} + \Lambda^{n-1T}\pd{G^{n-1}}{u^{n-2}} + \Lambda^{n-2T}\pd{G^{n-2}}{u^{n-2}} = 0, \\
...\\
\pd{J}{u^{1}}&= \pd{L}{u^{1}} + \Lambda^{2T}\pd{G^{2}}{u^{1}} + \Lambda^{1T}\pd{G^{1}}{u^{1}} = 0.\\
\end{aligned}
\end{equation}
Using $L = L(u^n,u^{n-1}, ...,u^{n-m}, D)$, we get the following equations:
\begin{equation}
\begin{aligned}
\pd{J}{u^{n}}  &= \pd{L}{u^{n}} + \Lambda^{nT}\pd{G^{n}}{u^{n}} = 0,\\
\pd{J}{u^{n-1}}&= \pd{L}{u^{n-1}} + \Lambda^{nT}\pd{G^{n}}{u^{n-1}} + \Lambda^{n-1T}\pd{G^{n-1}}{u^{n-1}} = 0,\\
...\\
\pd{J}{u^{n-m}}&=  \pd{L}{u^{n-m}} +\Lambda^{n-(m-1)T}\pd{G^{n-(m-1)}}{u^{n-m}} + \Lambda^{n-mT}\pd{G^{n-m}}{u^{n-m}} = 0,\\
\pd{J}{u^{n-(m+1)}}&=  \Lambda^{n-mT}\pd{G^{n-m}}{u^{n-(m+1)}} + \Lambda^{n-(m+1)T}\pd{G^{n-(m+1)}}{u^{n-(m+1)}} = 0,\\
...\\
\pd{J}{u^{1}}&= \Lambda^{2T}\pd{G^{2}}{u^{1}} + \Lambda^{1T}\pd{G^{1}}{u^{1}} = 0.\\
\end{aligned}
\end{equation}
Using the equation for the adjoint at the final pseudo-time step we obtain:
\begin{equation} \label{eq:AdjRFinal}
\left[\pd{G^n}{u^n}\right]^T\Lambda^n= -\left[\pd{L}{u^n}\right]^T, 
\end{equation}
using the other constraint equations and isolating the unknown adjoint term on the left hand side returns an adjoint recurrence relation for k = 2, 3, ...,m:
\begin{equation} \label{eq:AdjR_Source}
\pd{G^{k-1}}{u^{k-1}}^T \Lambda^{k-1}  = -\pd{G^{k}}{u^{k-1}}^T\Lambda^k.
\end{equation}
Similarly with the adjoint recurrence relation for k = m+1, ..., n-1 including a source term of the linearization of the objective with respect to the state at the given pseudo-time iteration, as follows:
\begin{equation} \label{eq:AdjR_Owin}
\pd{G^{k-1}}{u^{k-1}}^T \Lambda^{k-1}  = -\pd{G^{k}}{u^{k-1}}^T\Lambda^k -\left[\pd{L}{u^{k-1}}\right]^T.
\end{equation}
\\
For simulations where the objective function is only dependent on the final time-step --the focus of this work -- we only have one recurrence relation for all k = 2, 3, ..., n-1:
\begin{equation} \label{eq:AdjR_Gen}
\pd{G^{k-1}}{u^{k-1}}^T \Lambda^{k-1}  = -\pd{G^{k}}{u^{k-1}}^T\Lambda^k
\end{equation}
Lastly, we can take the derivative of equation (\ref{eq:AugObj}) with respect to the design variables to get the sensitivity equation.
\begin{equation} \label{eq:AdjLSens}
\td{J}{D} = \pd{L}{D} + \Lambda^{nT}\pd{G^{n}}{D} + \Lambda^{n-1T}\pd{G^{n-1}}{D} + \Lambda^{n-2T}\pd{G^{n-2}}{D} + ... 
\end{equation}

\subsection{Adjoint Computed Sensitivites for quasi-Newton Method}
For the quasi-Newton method we have the pseudo-time iteration given by equation (\ref{eq:PTCBDF1}), so that we define 
\begin{equation}
  G^k = G^k(u^k, u^{k-1}, D) = u^k - u^{k-1} + \left[P_{k-1} \right]^{-1}R(u^{k-1}),
\end{equation}
where $R = R(u^{k-1},x(D))$ denotes the residual as stated in section 2.2. We continue with the derivatives of constraint $G^k$, so that:
\begin{equation}
\begin{aligned}
\pd{G^k}{u^k} &= I, \\
\pd{G^k}{u^{k-1}} &= -I + \left[P_{k-1}\right]^{-1}\left[\pd{R(u^{k-1})}{u^{k-1}}\right]_2 + \pd{\left[P_{k-1}\right]^{-1}}{u^{k-1}}R(u^{k-1}), \\
\pd{G^k}{D} &= \left[P_{k-1}\right]^{-1}\pd{R(u^{k-1})}{D} + \pd{\left[P_{k-1}\right]^{-1}}{D}R(u^{k-1}).
\end{aligned}
\end{equation}
By using the differentiation of a matrix inverse shown in equation (\ref{eq:dKinv}) we can obtain
the derivative of the constraint term shown below: 
\begin{equation} \label{eq:Gcons}
\begin{aligned}
\pd{G^k}{u^k} &= I, \\
\pd{G^k}{u^{k-1}} &= -I + \left[P_{k-1}\right]^{-1}\left[\pd{R(u^{k-1})}{u^{k-1}}\right]_2 \\
&-\left[P_{k-1}\right]^{-1}\pd{\left[P_{k-1}\right]}{u^{k-1}} \left[P_{k-1}\right]^{-1}R(u^{k-1}),\\
\pd{G^k}{D} &= \left[P_{k-1}\right]^{-1}\pd{R(u^{k-1})}{D} - \left[P_{k-1}\right]^{-1}\pd{\left[P_{k-1}\right]}{D} \left[P_{k-1}\right]^{-1}R(u^{k-1}).\\ 
\end{aligned}
\end{equation}

Using the definition of the nonlinear solver increment in equation (\ref{eq:NM}) we can simplify the above
equations with:

\begin{equation} \label{eq:qNc}
\begin{aligned}
\pd{G^k}{u^k} &= I ,\\
\pd{G^k}{u^{k-1}} &= -I +
\left[P_{k-1}\right]^{-1}\left[\left[\pd{R(u^{k-1})}{u^{k-1}}\right]_2 - \pd{\left[P_{k-1}\right]}{u^{k-1}} \Delta u\right],\\
\pd{G^k}{D} &= \left[P_{k-1}\right]^{-1}\left[\pd{R(u^{k-1})}{D} - \pd{\left[P_{k-1}\right]}{x}\td{x}{D}\Delta u  \right].
\end{aligned}
\end{equation}
Please note that we can compute these Hessian vector products using complex
Frech\'et derivatives, rather than hand differentiating the residual operator
twice to obtain the Hessian operator; even though this is the adjoint linearization,
the hessian vector products are not transpose matrix vector products and can therefore be
computed using Frech\'et derivatives. The equation for the adjoint at the
final pseudo-time step  (\ref{eq:AdjRFinal}) with the first equation of the constraint derivatives (\ref{eq:Gcons}) gives the following initial source term:
\begin{equation} \label{eq:qNsAdj}
\left[I\right]\Lambda^n= -\left[\pd{L}{u^n}\right]^T.
\end{equation}
\\
Substituting in the constraint derivatives from equation (\ref{eq:qNc}) into equation (\ref{eq:AdjR_Gen}) returns:
\begin{equation}
\left[I\right]\Lambda^{k-1}  = -\left[-I + \left[P_{k-1}\right]^{-1}\left[\left[\pd{R(u^{k-1})}{u^{k-1}}\right]_2 - \pd{\left[P_{k-1}\right]}{u^{k-1}} \Delta u\right]   \right]^T\Lambda^k.
\end{equation}

We can write this recurrence relation in delta form as in the analysis solver:

\begin{equation} \label{eq:dNcAdjR}
\Delta \Lambda = -\left[\left[P_{k-1}\right]^{-1}\left[\left[\pd{R(u^{k-1})}{u^{k-1}}\right]_2 - \pd{\left[P_{k-1}\right]}{u^{k-1}} \Delta u\right]\right]^T\Lambda^k
\end{equation}
and distributing the transpose returns:
\begin{equation}
\Delta \Lambda = -\left[\left[\pd{R(u^{k-1})}{u^{k-1}}\right]_2 - \pd{\left[P_{k-1}\right]}{u^{k-1}} \Delta u \right]^T\left[P_{k-1} \right]^{-T} \Lambda^k.
\end{equation}
This motivates us to define a secondary adjoint variable for each recurrence
relation:
\begin{equation}
\left[P_{k-1} \right]^T\psi^k  = \Lambda^{k}.
\end{equation}
We then rewrite the delta form of the adjoint recurrence relation as follows:
\begin{equation} \label{eq:dPsiR}
\Delta \Lambda = -\left[\left[\pd{R(u^{k-1})}{u^{k-1}}\right]_2 - \pd{\left[P_{k-1}\right]}{u^{k-1}} \Delta u  \right]^T\psi^k.
\end{equation}
It is important to emphasize that, as previously noted, we do not need exact dual
correspondence here between the adjoint solver and the tangent one. Regardless of duality this formulation
will recover machine precision accuracy in the sensitivities in the limit of the linear system being solved to
machine precision. We substitute the constraint derivatives from equation (\ref{eq:qNc})
into the sensitivity equation (\ref{eq:AdjLSens}) and we obtain:
\begin{equation}
\begin{aligned}
\td{J}{D} &= \pd{L}{D} + \Lambda^{nT}\left[P_{n-1} \right]^{-1}\left[\pd{R(u^{n-1})}{D} - \pd{\left[P_{n-1}\right]}{x}\td{x}{D}\Delta u\right] \\
&+... \\
&+ \Lambda^{1T}\left[P_{0} \right]^{-1}\left[\pd{R(u^{0})}{D} - \pd{\left[P_{0}\right]}{x}\td{x}{D}\Delta
u\right].
\end{aligned}
\end{equation}
We can refer back to the definition of the secondary adjoint variable $\psi^k$ to simplify
this equation and remove an additional linear solve and get the equation below:

\begin{equation}
\begin{aligned}
\td{J}{D} &= \pd{L}{D} + \psi^{nT}\left[\pd{R(u^{n-1})}{D} - \pd{\left[P_{n-1}\right]}{x}\td{x}{D}\Delta u\right] \\
&+ ... \\
&+\psi^{1T}\left[\pd{R(u^{0})}{D} - \pd{\left[P_{0}\right]}{x}\td{x}{D}\Delta
u\right].
\end{aligned}
\end{equation}
We note that the adjoint sensitivity formulation uses only one approximate
linear solver per nonlinear step, like the forward and tangent solvers. To
guarantee convergence of the adjoint problem, we  use the dual solver of the
primal linear solver used at each nonlinear iteration. Please note that these
linearizations are only exact for cases in which the linear system is solved to
machine precision at each non-linear iteration. This is not a feasible
requirement for realistic CFD solvers, so we have to investigate the impact of
partial convergence of the linear system on the error in sensitivities.

\section{General Sensitivity Convergence Proof for Approximate Tangent Linearization of the Fixed Point Iteration}

We begin by defining the error in the objective function sensitivities as
\begin{equation} \label{eq:eps_l}
  \epsilon_L = \td{L}{D} - \widetilde{\td{L}{D}}, 
\end{equation}
where the tilde notation denotes an approximate value due to inexact linearization specifically. If we expand the above equation we recover
\begin{equation}
\begin{aligned}
  \epsilon^k_L &= \pd{L}{x}\td{x}{D} + \pd{L}{u^k}\td{u^k}{D} + \pd{L}{u^{k-1}}\td{u^{k-1}}{D} + ... + \pd{L}{u^{k-m}}\td{u^{k-m}}{D} \\
  &- \widetilde{\pd{L}{x}\td{x}{D}} - \widetilde{\pd{L}{u^k}\td{u^k}{D}} - \widetilde{\pd{L}{u^{k-1}}\td{u^{k-1}}{D}} - ... - \widetilde{\pd{L}{u^{k-m}}\td{u^{k-m}}{D}}, 
\end{aligned}
\end{equation}
which, when we factor in that all linearizations except for those of the fixed-point iteration are exact, reduces to the below equation
\begin{equation}
\begin{aligned}
  \epsilon^k_L &= \pd{L}{u^k}\left(\td{u^k}{D} - \widetilde{\td{u^k}{D}}\right) \\
  &+ \pd{L}{u^{k-1}}\left(\td{u^{k-1}}{D} - \widetilde{\td{u^{k-1}}{D}}\right)  \\
  &+ ... \\
  &+ \pd{L}{u^{k-m}}\left(\td{u^{k-m}}{D} - \widetilde{\td{u^{k-m}}{D}}\right). 
\end{aligned}
\end{equation}
For the purposes of increasing clarity in the tangent and adjoint proofs on the convergence of these inexactly linearized iterations, we proceed with all the proofs and numerical experiments in this work using objective functions computed only at the final iteration, even though the works that use the pseudo-time accurate adjoint do in fact use windowed objective functions. This simplifies $\epsilon^k_L$ into
\begin{equation} \label{eq:eps_lu}
  \epsilon^k_L = \pd{L}{u^k}\left(\td{u^k}{D} - \widetilde{\td{u^k}{D}}\right). 
\end{equation}
This motivates the definition of an expression for the error in the conservative variable sensitivities,
\begin{equation}
  \td{\epsilon^k_u}{D} = \td{u^k}{D} - \widetilde{\td{u^k}{D}},
\end{equation}
which simplifies (\ref{eq:eps_lu}) as follows
\begin{equation} \label{eq:eps_le}
  \epsilon^k_L = \pd{L}{u^k}\td{\epsilon^k_u}{D}. 
\end{equation}

We refer back to the definition of the fixed point iteration:
\begin{equation}
u^{k+1} = N(u^k, D) = u^k + H(u^k, D) = u^k + A(u^k,D)R(u^k,x(D))	
\end{equation}
where $A$ is some operator dependent on the pseudo-temporal discretization and $R$
is the appropriate residual operator for the spatial discretization of the
governing equations. For an exact linearization of the solution process we have
the expression below:
\begin{equation}
	\td{u^{k+1}}{D} = \td{N(u^{k})}{D} = \td{u^{k}}{D} +  \td{H}{D} = \td{u^{k}}{D} + \td{A}{D}R + A\td{R}{D}
\end{equation}
For an inexact linearization we inexactly linearize $A$ as $\widetilde{\td{A}{D}}$ and obtain:
\begin{equation}
		\widetilde{\td{u^{k+1}}{D}} = \widetilde{\td{N(u^{k})}{D}} = \widetilde{\td{u^{k}}{D}} + \widetilde{\td{A}{D}}R + A\widetilde{\td{R}{D}},
\end{equation}
To get the error in the conservative variable sensitivities we subtract the two expressions from one another:
\begin{equation}
	\td{u^{k+1}}{D} - \widetilde{\td{u^{k+1}}{D}} = \td{u^{k}}{D} - \widetilde{\td{u^{k}}{D}} + \td{A}{D}R - \widetilde{\td{A}{D}}R + A\td{R}{D} - A\widetilde{\td{R}{D}}.
\end{equation}
Using the definition of  $\td{\epsilon_u}{D}$ allows the grouping of like terms:
\begin{equation}
	\td{\epsilon_u^{k+1}}{D} = \td{\epsilon_u^{k}}{D} + \left[\td{A}{D} - \widetilde{\td{A}{D}}\right]R + A\left[\td{R}{D} - \widetilde{\td{R}{D}}\right].
\end{equation}
We then expand the derivative of the A matrix and the residual operator assuming that the partial derivatives of the residual operator are implemented correctly.
\begin{equation}
\begin{aligned}
	\td{\epsilon_u^{k+1}}{D} &= \td{\epsilon_u^{k}}{D} + \left[\pd{A}{x}\td{x}{D} + \pd{A}{u^k}\td{u^k}{D} - \widetilde{\pd{A}{x}\td{x}{D}} - \widetilde{\pd{A}{u^k}}\widetilde{\td{u^k}{D}}\right]R \\
	&+ A\left[\pd{R}{x}\td{x}{D} + \pd{R}{u^k}\td{u^k}{D} - \pd{R}{x}\td{x}{D} - \pd{R}{u^k}\widetilde{\td{u^k}{D}}\right].
\end{aligned}
\end{equation}

We can cancel out like terms to obtain:
\begin{equation}
\begin{aligned}
	\td{\epsilon_u^{k+1}}{D} &= \td{\epsilon_u^{k}}{D} + \left[\pd{A}{x}\td{x}{D} + \pd{A}{u^k}\td{u^k}{D} - \widetilde{\pd{A}{x}}\td{x}{D} - \widetilde{\pd{A}{u^k}}\widetilde{\td{u^k}{D}}\right]R \\
	&+ A\left[\pd{R}{u^k}\td{u^k}{D} - \pd{R}{u^k}\widetilde{\td{u^k}{D}}\right]
\end{aligned}
\end{equation}

We can then group terms and use the definition of $\td{\epsilon_u^k}{D}$:
\begin{equation}
\begin{aligned}
  \td{\epsilon_u^{k+1}}{D} &= \td{\epsilon_u^{k}}{D} + \left[\left(\pd{A}{x} - \widetilde{\pd{A}{x}}\right)\td{x}{D} + \left( \pd{A}{u^k} - \widetilde{\pd{A}{u^k}}\right)\td{u^k}{D}  + \pd{A}{u^k}\td{\epsilon_u^k}{D} - \left(\pd{A}{u^k} - \widetilde{\pd{A}{u^k}} \right) \td{\epsilon_u^k}{D} \right]R \\
  &+ A\left[\pd{R}{u^k}\td{\epsilon_u^k}{D}\right],
\end{aligned}
\end{equation}

we can then rearrange terms:

\begin{equation} \label{eq:TanConsErr}
\begin{aligned}
	\td{\epsilon_u^{k+1}}{D} &= \td{\epsilon_u^{k}}{D} + A\left[\pd{R}{u^k}\td{\epsilon_u^k}{D}\right] + \pd{A}{u^k}\td{\epsilon_u^k}{D}R \\
  &+ \left[\left(\pd{A}{x} - \widetilde{\pd{A}{x}}\right)\td{x}{D} + \left( \pd{A}{u^k} - \widetilde{\pd{A}{u^k}}\right)\td{u^k}{D}  - \left(\pd{A}{u^k} - \widetilde{\pd{A}{u^k}} \right)\td{\epsilon_u^k}{D} \right]R, \\
\end{aligned}
\end{equation}

and we can use triangle inequality to obtain the below, where the norm used in this work is the 2-norm: 

\begin{equation} \label{eq:TriangIneq}
\begin{aligned}
    \norm{\td{\epsilon_u^{k+1}}{D}} &<   \norm{\td{\epsilon_u^{k}}{D} + A\left[\pd{R}{u^k}\td{\epsilon_u^k}{D}\right] + \pd{A}{u^k}\td{\epsilon_u^k}{D}R} \\
    &+ \norm{\left[\left(\td{A}{D} - \widetilde{\td{A}{D}}\right) - \left(\pd{A}{u^k} - \widetilde{\pd{A}{u^k}} \right)\td{\epsilon_u^k}{D}   \right]R} .
\end{aligned}
\end{equation}

We can group the first three terms on the right hand side into $B$ where $B = \pd{N}{u}$, by Cauchy-Schwarz inequality we know that:
\begin{equation} \label{eq:CS_ineq}
\norm{\td{\epsilon_u^{k}}{D} + A\left[\pd{R}{u^k}\td{\epsilon_u^k}{D}\right] + \pd{A}{u^k}\td{\epsilon_u^k}{D}R} < \norm{B}\norm{\td{\epsilon_u^k}{D}}
\end{equation}

$B$ is the derivative of the contractive fixed point iteration, therefore
$\norm{B} < 1$. Using the inequality (\ref{eq:CS_ineq}) in the inequality (\ref{eq:TriangIneq}) returns:  

\begin{equation}
\begin{aligned}
    \norm{\td{\epsilon_u^{k+1}}{D}} < \norm{B}\norm{\td{\epsilon_u^k}{D}}  
    + \norm{\left[\left(\td{A}{D} - \widetilde{\td{A}{D}}\right)  - \left(\pd{A}{u^k} - \widetilde{\pd{A}{u^k}} \right)\td{\epsilon_u^k}{D}\right]R}.
\end{aligned}
\end{equation}
As a result, the error in the conservative variable sensitivities expressed by
$\td{\epsilon_u^k}{D}$ decreases as the residual decreases and the contractivity of
the fixed point iteration progresses through the primal solution process. The triangle inequality is used to show the pseudo-temporal evolution of the error as governed by the inequality below.
We first refer back to equation (\ref{eq:eps_lu}) to get
\begin{equation} \label{eq:TanErrSensIneq}
  \norm{\epsilon^{k+1}_L} < \norm{\pd{L}{u^{k+1}}}\norm{\td{\epsilon_u^{k+1}}{D}},
\end{equation}
and so we obtain the final inequality to define $\epsilon^k_L$, or the error in the objective function sensitivities:
\begin{proposition} \label{prop:TanErrIneq}
\begin{equation}
    \norm{\epsilon^{k+1}_L} < \norm{\pd{L}{u^{k+1}}}\left[\norm{B}\norm{\td{\epsilon^k}{D}} 
    + \norm{\left[\left(\td{A}{D} - \widetilde{\td{A}{D}}\right)  - \left(\pd{A}{u^k} - \widetilde{\pd{A}{u^k}} \right)\td{\epsilon_u^k}{D}  \right]R}\right],
\end{equation}

where we have regrouped the last term of the right hand side of equation (\ref
{eq:TanConsErr}). This shows that once the residual is much lower than the
sensitivity error, the sensitivity converges as the contractivity of the
fixed-point iteration of the nonlinear problem. This does contain a lagging
effect that has been previously observed in work on the one-shot Adjoint method
in "piggy-back iterations" \cite{GR_OSA}. One can also note that in cases
where $A$ approaches $-\pd{R}{u}^{-1}$ there is no dependence on contractivity of
the fixed point in the derivative and the error will be multiplied by the
residual in all terms and will converge as the nonlinear problem does. 
\end{proposition}

\begin{corollary}
  The smaller the error in the linearization of $A$, signified by $\td{A}{D} - \widetilde{\td{A}{D}}$, the earlier in iteration-space the contractivity of the fixed-point dominates the convergence of the sensitivities. Note that in Newton-type solvers using the inverse identity in equation (\ref{eq:dKinv}) the error is directly a function of the linear system tolerance.
\end{corollary}
\begin{corollary}
Should we desire to analyze the error convergence of the sensitivities of a windowed objective function, it would entail the addition of:
 \begin{equation}
 \sum^m_{j=1} \norm{\pd{L}{u^{k-j}}}\norm{\td{\epsilon_u^{k-j}}{D}}
 \end{equation}
 into the right side of equation \ref{eq:TanErrSensIneq} in a straightforward manner.
\end{corollary}
\section{General Sensitivity Convergence Proof for Approximate Adjoint Linearization of the Fixed Point Iteration}
As in the tangent section, we wish to quantify the error introduced to the
sensitivity calculation by the approximate linearization of the fixed point
iteration. We know that for a properly implemented linearization and
transposition that the adjoint and tangent converge to the same sensitivity values
(even for approximate linearizations) and therefore that the error in the
sensitivities goes to zero as the nonlinear problem is converged. Beginning from
equation (\ref{eq:AdjLSens}) first shown in the initial pseudo-time accurate adjoint
derivations and reproduced below, where $G^k$ is the fixed point iteration shifted by its output,
$G^k = u^k - N(u^{k-1}, x(D))$, below:

\begin{equation}
\td{L}{D} = \pd{L}{D} + \Lambda^{nT}\pd{G^{n}}{D} + \Lambda^{n-1T}\pd{G^{n-1}}{D} + \Lambda^{n-2T}\pd{G^{n-2}}{D} + ... + \Lambda^{1T}\pd{G^1}{D}.
\end{equation}

As in the tangent proof we have the following identity for the fixed point iteration. It is important to note that by definition of the fixed point, $A$ is not orthogonal to R and it is bounded away from 0, lest the fixed point terminate at a state that does not satisfy the discretized form of the governing equations (signifed by $R = 0$).
\begin{equation}
	u^{k+1} = N(u^k, D) = u^k + H^{k+1}(u^k, x(D)) = u^k + A(u^k,x(D))R(u^k,x(D)).
\end{equation}

The fixed point iteration derivatives are expressed below, where the linearization of the residual is assumed to be exact, but the linearization of the $A$ matrix is approximate in the actual implementation:
\begin{equation}
\begin{aligned}
	\pd{H^{k+1}}{u^k} & = \pd{A(u^k, x(D))}{u^k}R(u^k, x(D)) + A(u^k, x(D))\pd{R(u^k, x(D))}{u^k}, \\
	\pd{u^{k+1}}{D} &` = \pd{A(u^k, x(D))}{D}R(u^k, x(D)) + A(u^k, x(D))\pd{R(u^k, x(D))}{D}.
\end{aligned}
\end{equation}

Substituting in the approximate terms into (\ref{eq:AdjLSens}) we get the expression below:

\begin{equation}
\widetilde{\td{L}{D}} = \pd{L}{D} + \tilde{\Lambda}^{nT}\widetilde{\pd{G^{n}}{D}} + \tilde{\Lambda}^{n-1T}\widetilde{\pd{G^{n-1}}{D}} + \tilde{\Lambda}^{n-2T}\widetilde{\pd{G^{n-2}}{D}} + ... + \tilde{\Lambda}^{1T} \widetilde{\pd{G^{1}}{D}} .
\end{equation}

We can subtract the exact sensitivity equation from the approximate one:

\begin{equation}
\begin{aligned}
	\td{L}{D} - \widetilde{\td{L}{D}} &= \left(\pd{L}{D} - \pd{L}{D}\right) \\
	&+ \left(\Lambda^{nT}\pd{G^{n}}{D} - \tilde{\Lambda}^{nT}\widetilde{\pd{G^{n}}{D}}\right) \\
	&+ \left(\Lambda^{n-1T}\pd{G^{n-1}}{D} - \tilde{\Lambda}^{n-1T}\widetilde{\pd{G^{n-1}}{D}}\right) \\
	&+ \left(\Lambda^{n-2T}\pd{G^{n-2}}{D} - \tilde{\Lambda}^{n-2T}\widetilde{\pd{G^{n-2}}{D}}\right) \\
    &+ ... \\
   	&+ \left(\Lambda^{1T}\pd{G^{1}}{D} - \tilde{\Lambda}^{1T}\widetilde{\pd{G^{1}}{D}}\right) . \\
\end{aligned}
\end{equation}

We define three error terms for the errors at each nonlinear iteration in the
pseudo-time adjoint, the linearization of the nonlinear iteration with respect
to the state variable and the linearization with respect to the design
variables, denoted by $\epsilon^k_\Lambda, \epsilon^k_u, \epsilon^k_D$ respectively, and given by:

\begin{equation} \label{eq:epsdeforig}
\begin{aligned}
	\epsilon^k_\Lambda &= \Lambda^k - \tilde{\Lambda}^k, \\
	\epsilon^k_u &= \pd{G^k}{u^{k-1}} - \widetilde{\pd{G^k}{u^{k-1}}}, \\
	\epsilon^k_D &= \pd{G^k}{D} - \widetilde{\pd{G^k}{D}}.
\end{aligned}
\end{equation}

Using these epsilon definitions we can simplify the equations into the below equation with the epsilon terms as the unknowns:

\begin{equation} \label{eq:SensErr}
\begin{aligned}
	\td{L}{D} - \widetilde{\td{L}{D}} &= \epsilon^{nT}_\Lambda\pd{G^n}{D} + \Lambda^{nT}\epsilon^n_D - \epsilon^{nT}_\Lambda\epsilon_D^n\\
	&+ \epsilon^{n-1T}_\Lambda\pd{G^{n-1}}{D} + \Lambda^{n-1T}\epsilon^{n-1}_D - \epsilon^{n-1T}_\Lambda\epsilon_D^{n-1} \\
	&+ \epsilon^{n-2T}_\Lambda\pd{G^{n-2}}{D} + \Lambda^{n-2T}\epsilon^{n-2}_D - \epsilon^{n-2T}_\Lambda\epsilon_D^{n-2} \\
    &+ ... \\
   	&+ \epsilon^{1T}_\Lambda\pd{G^{1}}{D} + \Lambda^{1T}\epsilon^{1}_D - \epsilon^{1T}_\Lambda\epsilon_D^{1}.\\
\end{aligned}
\end{equation}
which can be rewritten using summation notation and the definition of $\epsilon^k_L$ in equation (\ref{eq:eps_l})
\begin{equation} \label{eq:nkSensErr}
\begin{aligned}
  \epsilon^n_L &= \sum_{k=1}^n\epsilon^{kT}_\Lambda\pd{G^k}{D} + \Lambda^{kT}\epsilon^k_D - \epsilon^{kT}_\Lambda\epsilon_D^k \\
  &= \sum_{k=1}^n \Lambda^{kT}\epsilon^k_D + \epsilon^{kT}_\Lambda\widetilde{\pd{G}{D}}^k
\end{aligned}
\end{equation}
It is clear that for the error analysis to parallel that of the tangent section that $\epsilon^{k}_\Lambda$ and $\epsilon^k_D$ must converge at the same rate as the convergence of the primal problem. As such, this proof cannot proceed any further without expressions for the epsilon error terms, and we proceed using the fact that only the $\pd{A}{()}$ term has an approximation to get the following two identities:

\begin{equation} \label{eq:eps_def}
\begin{aligned}
	\epsilon^k_u &= \pd{G^k}{u^{k-1}} - \widetilde{\pd{G^k}{u^{k-1}}} = \pd{A^k}{u^{k-1}}R + A\pd{R}{u^{k-1}} - \widetilde{\pd{A^k}{u^{k-1}}}R - A\pd{R}{u^{k-1}} \\
  &= \left(\pd{A^k}{u^{k-1}} - \widetilde{\pd{A^k}{u^{k-1}}}   \right)R, \\
	\epsilon^k_D &= \pd{G^k}{D} - \widetilde{\pd{G^k}{D}} = \pd{A}{D}R + A\pd{R}{D} - \widetilde{\pd{A}{D}}R - A\pd{R}{D} = \left(\pd{A}{D} - \widetilde{\pd{A}{D}}   \right)R.
\end{aligned}
\end{equation}

The two error terms above have a scaling with the residual that is important to
keep in mind. We use the definitions of the error terms to get an expression of
the difference between the sensitivities with exact and approximate
linearizations. Using equation (\ref{eq:AdjR_Gen}) gives an error recurrence
relationship:
\begin{equation} \label{eq:erradj_rec}
\begin{aligned}
	\epsilon^{k-1}_\Lambda &= -\pd{G^{k}}{u^{k-1}}^T\Lambda^k + \widetilde{\pd{G^{k}}{u^{k-1}}^T}\tilde{\Lambda}^k \\
  &=-\epsilon_u^{kT}\Lambda^k - \widetilde{\pd{G^k}{u^{k-1}}^T}\epsilon^k_{\Lambda} .
\end{aligned}
\end{equation}

Applying this equation to the second to last iteration we have

\begin{equation}
  \epsilon^{n-1}_\Lambda = -\epsilon_u^{nT}\Lambda^n + \epsilon_u^{nT}\epsilon_\Lambda^n - \pd{G^n}{u^{n-1}}\epsilon_\Lambda^n . 
\end{equation}

Since $\Lambda^n = \pd{L}{u^n}$ which has no approximation error, $\epsilon_\Lambda^n = 0$. Substituting in the expression for $\epsilon_u^{kT}$ we get:
\begin{equation}
  \epsilon^{n-1}_\Lambda = -\left(\pd{A}{u} - \widetilde{\pd{A}{u}}   \right)R^n\Lambda^n .
\end{equation}

Should we proceed through the same steps for $\epsilon^{n-2}_\Lambda$, we obtain:
\begin{equation}
  \epsilon^{n-2}_\Lambda = -\left(\pd{A}{u} - \widetilde{\pd{A}{u}}   \right)R^{n-1}\Lambda^{n-1} .
\end{equation}

It is therefore clear that in order to show that the error goes to zero it must be
shown that $\norm{R^k\Lambda^k} \approx 0$ for all iterations. It is important
here to refer back to the definition of the fixed point iteration; A is bounded
away from 0 and A is not orthogonal to R (lest the fixed point terminate at a
state that does not satisfy the discretized form of the PDE). Therefore
$\norm{A}\neq0$ and it is possible to move to an analysis of the $\Lambda^kR^k$
term. Since the pseudo-temporal adjoint converges at the reverse of the analysis
process (through transpose and linearizing the fixed point iteration),
$\norm{H^k\Lambda^k} = \norm{H^n\Lambda^n}$, which for a converged simulation is
on the order of machine zero. Since A is bounded away from zero and not
orthogonal to R, this means that $\norm{R^k\Lambda^k} \approx 0$, and all the
errors in the linearization of the A operator (denoted by $\epsilon_u$) do not
contribute to the sensitivity error. The reverse convergence of the adjoint as
compared to the analysis mode is confirmed by the results in \cite{PTA-iqn},
where the adjoint is shown to converge to its final value in the reverse of the
analysis problem. Thus 

\begin{equation}
	\epsilon^{k-1}_\Lambda = -\left(\pd{A}{u} - \widetilde{\pd{A}{u}}   \right)R^k\Lambda^k \approx 0 .
\end{equation}

Otherwise, the error at every iteration is scaled by a magnitude of the final state residual for a stalled or truncated simulation. Using that $\epsilon_\Lambda^k \approx 0$ for a converged simulation we get the following expression when neglecting the terms in \ref{eq:SensErr} that are multiplied by $\Lambda^kR^k \approx 0$:

\begin{equation}
	\td{L}{D} - \widetilde{\td{L}{D}} = \Lambda^{nT}\epsilon^n_D + \Lambda^{n-1T}\epsilon^{n-1}_D + \Lambda^{n-2T}\epsilon^{n-2}_D + ... + \Lambda^{1T}\epsilon^{1}_D .
\end{equation}

Substituting in the expression for the error in the linearization of the nonlinear iteration with respect to the design variable returns the following:
\begin{equation}
\begin{aligned}
	\td{L}{D} - \widetilde{\td{L}{D}} &= -\Lambda^{nT}\left(\pd{A}{u}^n - \widetilde{\pd{A}{u}}^n   \right)R^n \\
	&- \Lambda^{n-1T}\left(\pd{A}{u}^{n-1} - \widetilde{\pd{A}{u}}^{n-1}   \right)R^{n-1} \\
	&- \Lambda^{n-2T}\left(\pd{A}{u}^{n-2} - \widetilde{\pd{A}{u}}^{n-2}   \right)R^{n-2} \\
	&- ... \\
	&- \Lambda^{1T}\left(\pd{A}{u}^{1} - \widetilde{\pd{A}{u}}^{1}   \right)R^{1} .
\end{aligned}
\end{equation}

\begin{proposition} \label{prop:AdjErrEq}

We can rearrange the above expression and use the same argument as above to show that $\norm{R^k\Lambda^k} \approx \norm{\Lambda^nR^n} \approx 0$ to obtain the final expression for the error:

\begin{equation}
\begin{aligned}
  \td{L}{D} - \widetilde{\td{L}{D}} &= -\left(\pd{A}{u}^n - \widetilde{\pd{A}{u}}^n   \right)^T\Lambda^{n}R^n \\
  &- \left(\pd{A}{u}^{n-1} - \widetilde{\pd{A}{u}}^{n-1}   \right)^T\Lambda^{n-1}R^{n-1} \\
  &- \left(\pd{A}{u}^{n-2} - \widetilde{\pd{A}{u}}^{n-2}   \right)^T\Lambda^{n-2}R^{n-2} \\
  &- ... \\
  &- \left(\pd{A}{u}^{1} - \widetilde{\pd{A}{u}}^{1}   \right)^T\Lambda^{1}R^{1} \approx 0.
\end{aligned}
\end{equation}
\end{proposition}

Again we see that the convergence of the nonlinear problem leads to less error
in the sensitivities even in unconverged flows. From the expressions above we
can see that the error in the sensitivities is dependent on the error in the A
operator linearization and the residual of the nonlinear problem. 
\begin{corollary} 
 For quasi-Newton solvers like the ones that are the primary
 focus of this work, this indicates that the convergence is a multiple of the
 tolerance of the linear system -- as the error in the linearization of the
 matrix inverse scales with the linear system tolerance -- and the convergence
 of the non-linear problem, and we will see this behavior borne out in the
 results section.
\end{corollary}
\section{Impact of Appromixate Linearization on Sensitivity Accuracy}
\begin{figure}[!h]
  \centering
  \includegraphics[keepaspectratio, width=.75\textwidth]{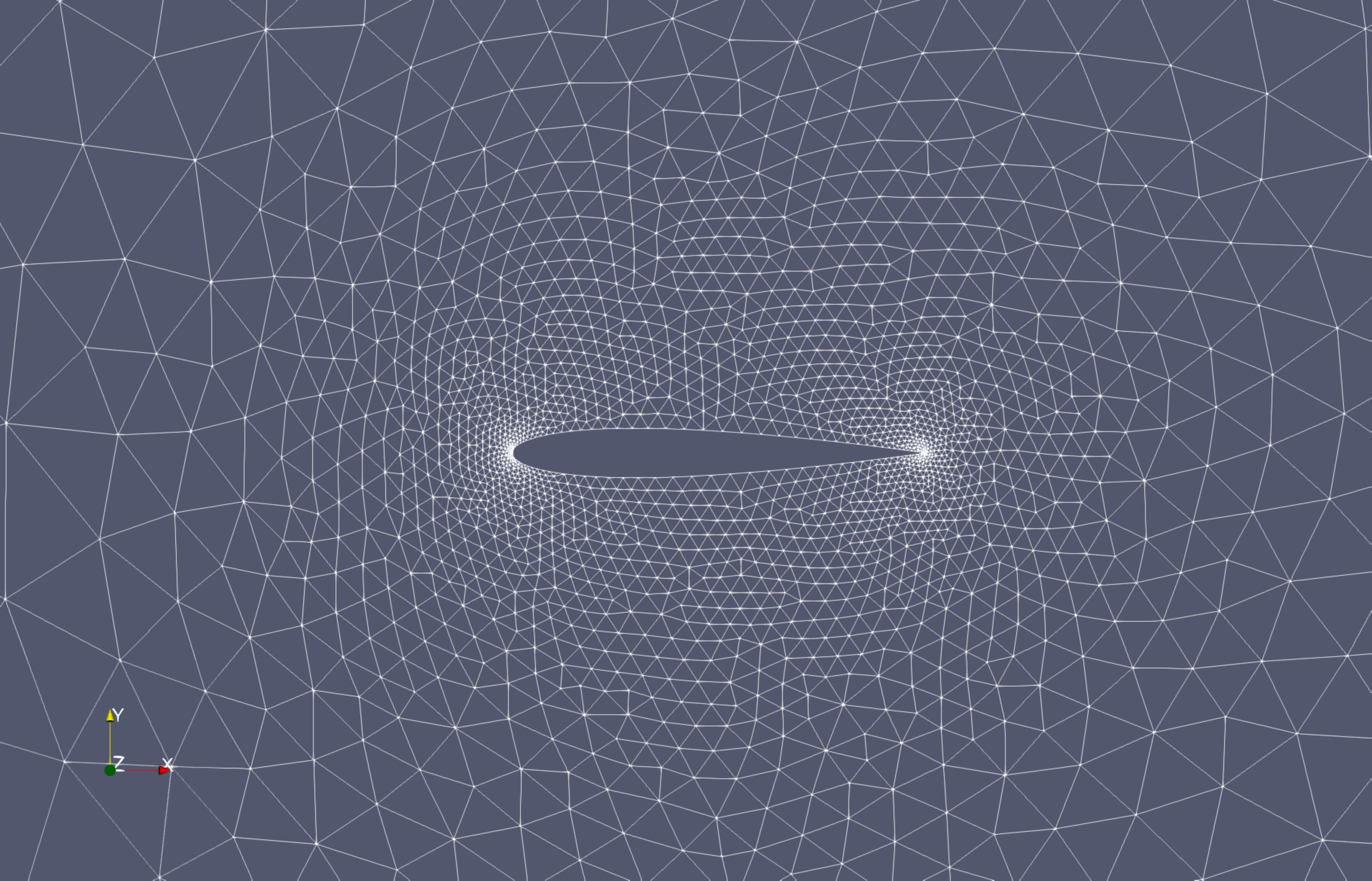}
\caption{Computational mesh for NACA0012 airfoil \label{fig:mesh}.}
\end{figure}

This section shows the effect of an approximate linearization on the accuracy of
the sensitivity computation. The first two sets of results portray the effect of
partial linear solves on the accuracy  of the identity of the differentiation of
the inverse matrix in the sensitivity computation shown in equation
(\ref{eq:dKinv}). The final set shows the impact of exactly linearizing the
residual operator at every stage of a Runge-Kutta scheme when the gradients are only
computed at the first stage and then frozen throughout the analysis.  All cases
were run on an unstructured triangular mesh consisting of 4212 elements shown in
Figure (\ref{fig:mesh}), in $M=.7$ flow with $\alpha = 2^o$, with 2 Hickes-Henne bump functions \cite{HicksHenneDesign} used as design variables to perturb the airfoil surface.

\subsection{Results for an Exact Jacobian Augmented with a Mass Matrix}
Here we look at the impact of inexact linearization for an exact Newton
algorithm, when the left hand side is the exact linearization of the residual
operator augmented with a suitable mass matrix. The fixed point iteration is
shown below,

\begin{equation}
\begin{aligned}
u^{k} &= u^{k-1} + \Delta u ,
\end{aligned}
\end{equation} 
where $\Delta u$ is computed by solving the linear system below to a linear tolerance of the user's choice; in this work  FGMRES algorithm \cite{GMRES,PTA-iqn} is used to solve said linear system, except for the cases where we show duality, where we use Gauss-Seidel sweeps as they are right-hand-side independent.
\begin{equation}
\left[P_{k-1}\right]\Delta u = -R ,
\end{equation}
where R is the residual operator of the discretization and $\left[P_{k-1}\right]$ is a left hand side matrix that determines the solver. As this is an exact Newton algorithm, $\left[P_{k-1}\right]$ is the Jacobian of
the spatial discretization augmented with a suitable mass matrix, shown below and defined in the background section, given by
\begin{equation} 
\left[P_{k-1}\right] = \left[\pd{R}{u^{k-1}}\right] + \frac{vol}{\Delta t CFL}.
\end{equation}

\begin{figure}
\begin{subfigmatrix}{2}
  \subfigure[Design Variable
1]{\includegraphics[keepaspectratio]{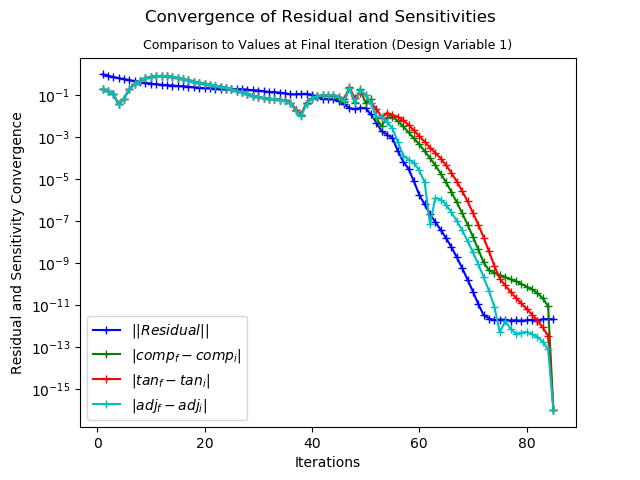}}
  \subfigure[Design Variable
2]{\includegraphics[keepaspectratio]{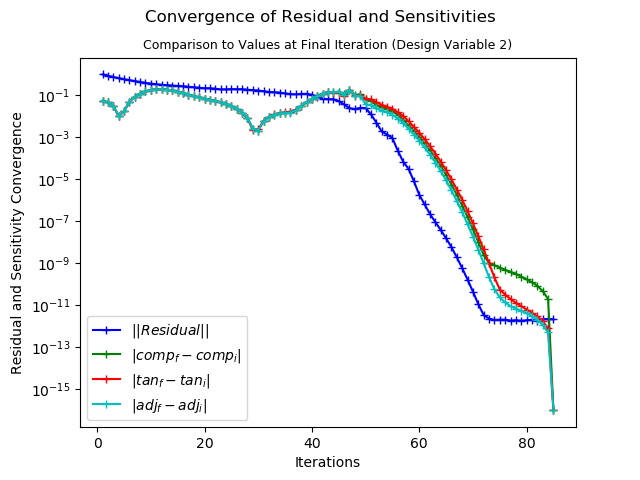}}
\end{subfigmatrix}
\caption{Convergence of nonlinear problem and sensitivities for linear tolerance of $1e-1$:
measured by $L_2$-norm of the residual and the difference between current and final sensitivity values.}
\label{fig:SensConv_LT1e-1}
\end{figure}

\begin{figure}
\begin{subfigmatrix}{2}
  \subfigure[Design Variable
1]{\includegraphics[keepaspectratio]{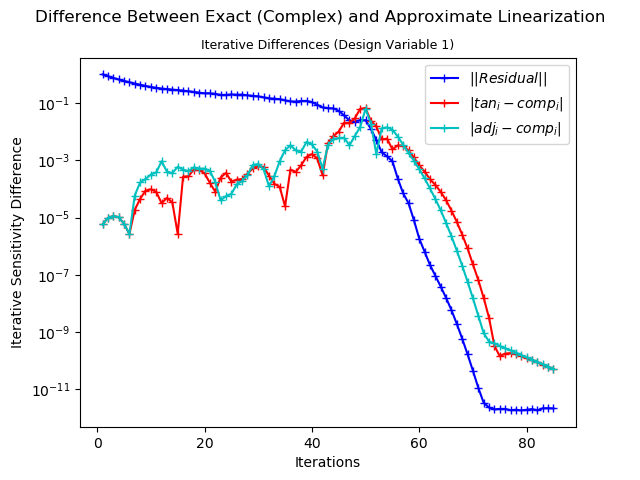}}
  \subfigure[Design Variable
2]{\includegraphics[keepaspectratio]{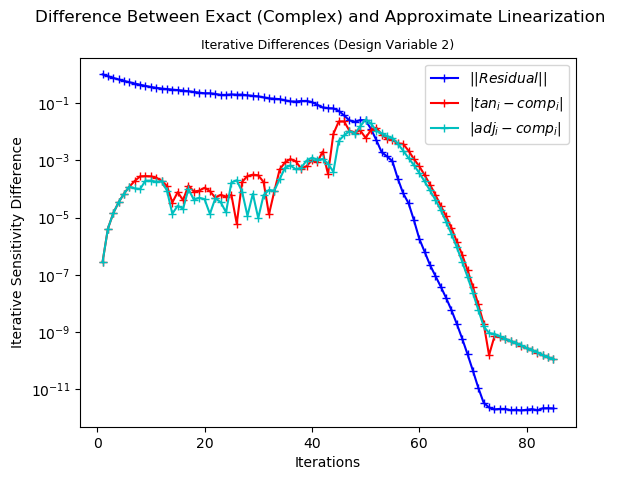}}
\end{subfigmatrix}
\caption{Difference between adjoint, tangent, and complex sensitivities at each iteration for a linear tolerance of $1e-1$.}
\label{fig:SensComp_LT1e-1}
\end{figure}

We can see in Figure (\ref{fig:SensConv_LT1e-1}) that for the linear tolerance
$1e-1$  the complex, tangent, and adjoint sensitivities, converge to their final values
at the same rate as the analysis problem itself, which is expected based on the
formulation. We can then see in Figure (\ref{fig:SensComp_LT1e-1}), which
depicts the difference between the complex and tangent sensitivities
over the iteration history of the analysis solution process, that the maximum
difference is of the order of the linear tolerance of the linear system.
Furthermore, as the analysis problem converges, so do the complex, tangent, and adjoint
sensitivities to each other despite the inexact differentiation. We can see that
at full convergence of the analysis problem, the adjoint, tangent and complex-step
sensitivities correspond to each other to a high degree of precision and these
would also correspond to the steady-state tangent and adjoint computed
sensitivities linearized about the converged analysis state.

\begin{figure}
\begin{subfigmatrix}{2}
  \subfigure[Design Variable
1]{\includegraphics[keepaspectratio]{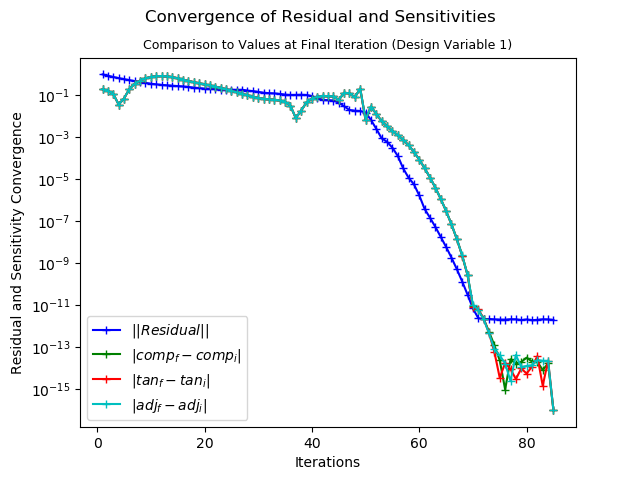}}
  \subfigure[Design Variable
2]{\includegraphics[keepaspectratio]{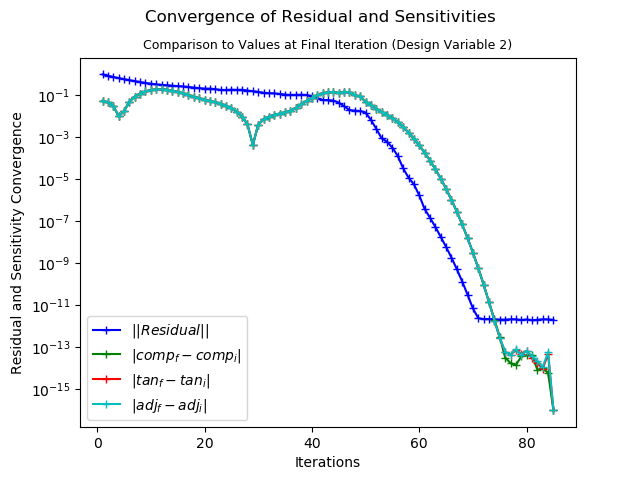}}
\end{subfigmatrix}
\caption{Convergence of nonlinear problem and sensitivities for linear tolerance of $1e-4$:
measured by $L_2$-norm of the residual and the difference between current and final sensitivity values.}
\label{fig:SensConv_LT1e-4}
\end{figure}

\begin{figure}
\begin{subfigmatrix}{2}
  \subfigure[Design Variable
1]{\includegraphics[keepaspectratio]{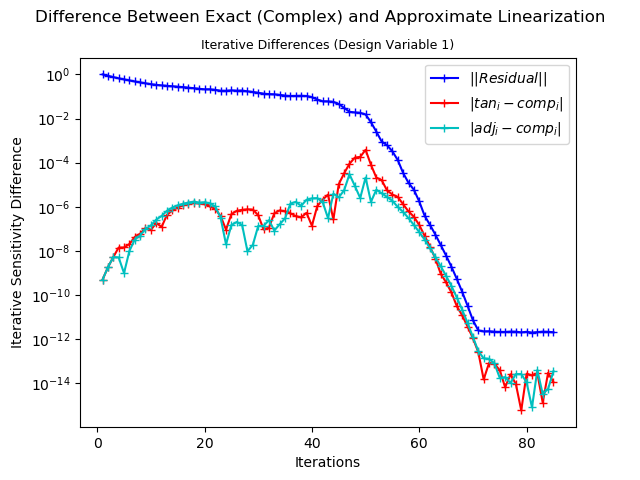}}
  \subfigure[Design Variable
2]{\includegraphics[keepaspectratio]{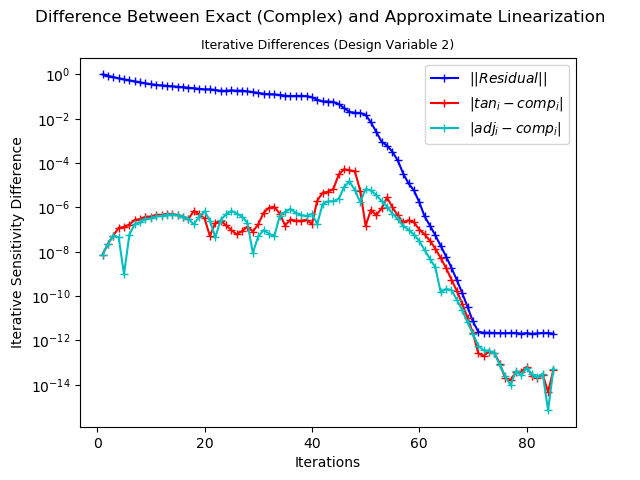}}
\end{subfigmatrix}
\caption{Difference between adjoint, tangent, and complex sensitivities at each iteration for a linear tolerance of $1e-4$.}
\label{fig:SensComp_LT1e-4}
\end{figure}

Figure (\ref{fig:SensConv_LT1e-4}) contains the same information as Figure
(\ref{fig:SensConv_LT1e-1}), and Figure (\ref{fig:SensComp_LT1e-1}) is the
sister plot of Figure (\ref{fig:SensComp_LT1e-4}) but the two later plots show
results with a tighter linear system tolerance of $1e-4$. We can see that as
before the maximum iterative difference is again on the order of the linear
system tolerance, and that as the analysis problem converges so do the adjoint, tangent,
and complex sensitivities to each other, down to nearly machine precision. 

One thing to note is that while the adjoint and tangent linearizations have very
similar convergence behavior, they are not identical. That is because these
simulations use a specified linear tolerance to terminate the linear solver,
and the number of steps is not identical between the tangent and adjoint
linearizations. The goal of this investigation is to show similar behavior
despite the loss of strict duality as identical behavior has already been shown
for cases when strict duality is maintained \cite{PTA-devel}.

\begin{figure}
 \begin{subfigmatrix}{2}
    \centering
    \subfigure[Maximum Iterative Difference]{\includegraphics[keepaspectratio, width=.45\textwidth]{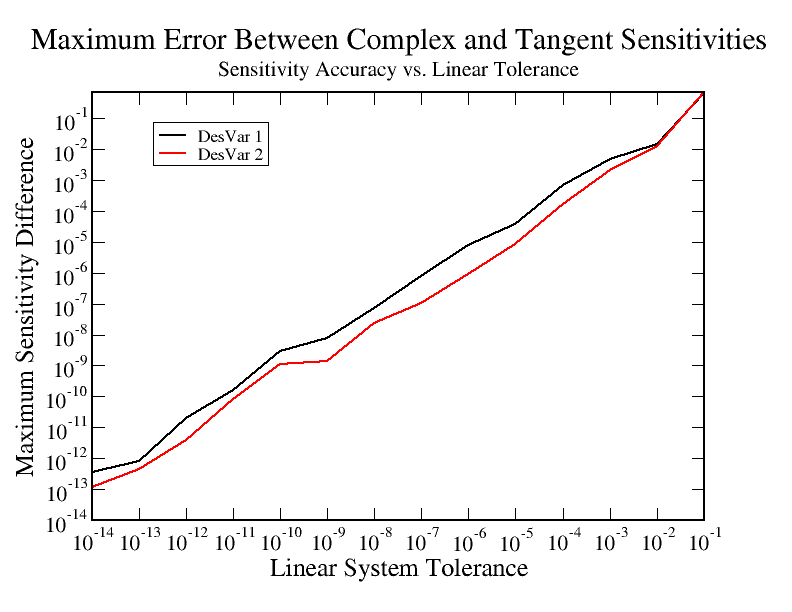}}
    \subfigure[Minimum Iterative Difference]{\includegraphics[keepaspectratio, width=.45\textwidth]{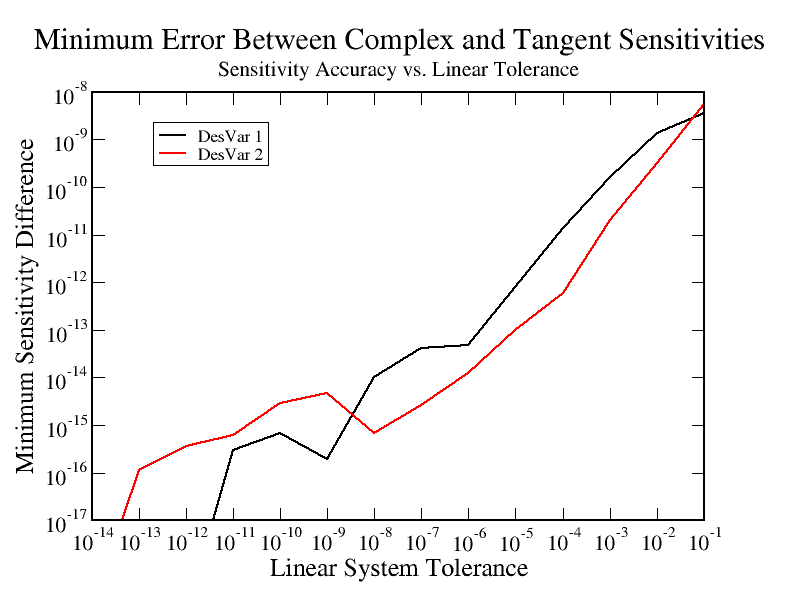}}
\end{subfigmatrix}   
\caption{Maximum and minimum difference between tangent and complex sensitivities as a function of linear tolerance}.
\label{fig:ErrConv_QN}
\end{figure}

Having seen the impact of the tighter linear system tolerance on the maximum
iterative difference between the tangent and complex-step finite-difference computed sensitivities, we simulate
the analysis problem with linear tolerances at every order from $1e-1$ to $1e-14$, and
plot the maximum iterative difference in Figure (\ref{fig:ErrConv_QN}). The
maximum iterative difference is directly related to the linear system tolerance.
This allows for good estimates of the maximum iterative error as a function of
the linear system tolerance. The minimum iterative difference shows similar
behavior, but it is a function of the linear tolerance and the convergence of
the non-linear problem.

\FloatBarrier

\subsection{Results for an Inexact Jacobian Augmented with a Mass Matrix}
Here we show the impact of inexact linearization for an inexact-quasi-Newton
algorithm, i.e. the left hand side is the first order linearization of the
second order residual operator augmented with a suitable mass matrix. Similar to the previous section the nonlinear solver is defined by

\begin{equation} \label{eq:PTCBDF}
\begin{aligned}
u^{k} &= u^{k-1} + \Delta u
\end{aligned}
\end{equation} 
where $\Delta u$ is computed by solving the linear system 
\begin{equation}
\left[P_{k-1}\right]\Delta u = -R
\end{equation}
to a linear tolerance of the user's choice using the same algorithms as in the exact Newton section.

However, here $\left[P_{k-1}\right]$ is an inexact linearization; it is the Jacobian of the first order spatial discretization
augmented with a diagonal term to ensure that it is diagonally dominant, shown
in equation (\ref{eq:PC2}):
\begin{equation} \label{eq:PC2}
\left[P_{k-1}\right] = \left[\pd{R}{u^{k-1}}\right]_1 + \frac{vol}{\Delta t CFL}.
\end{equation}

\begin{figure}
\begin{subfigmatrix}{2}
  \subfigure[Design Variable
1]{\includegraphics[keepaspectratio]{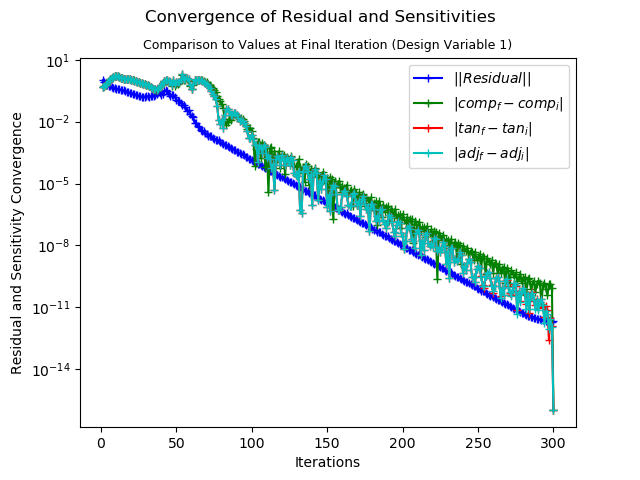}}
  \subfigure[Design Variable
2]{\includegraphics[keepaspectratio]{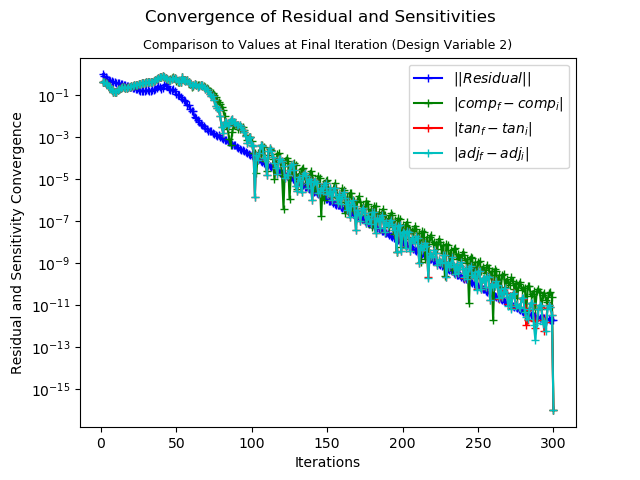}}
\end{subfigmatrix}
\caption{Convergence of nonlinear problem and sensitivities for linear tolerance of $1e-1$:
measured by $L_2$-norm of the residual and the difference between current and final sensitivity values}
\label{fig:2ndSensConv_LT1e-1}
\end{figure}

\begin{figure}
\begin{subfigmatrix}{2}
  \subfigure[Design Variable
1]{\includegraphics[keepaspectratio]{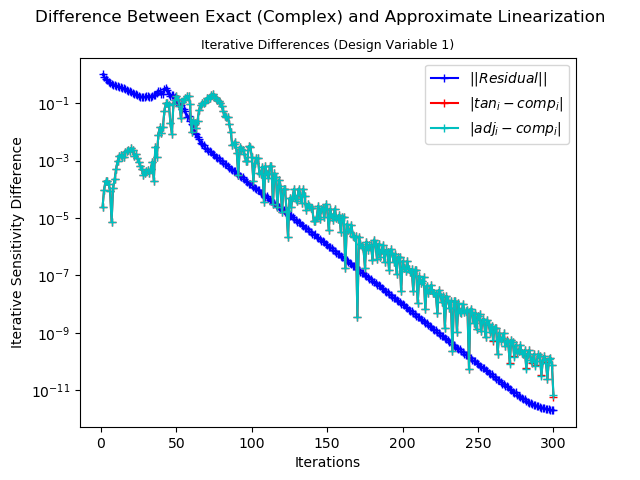}}
  \subfigure[Design Variable
2]{\includegraphics[keepaspectratio]{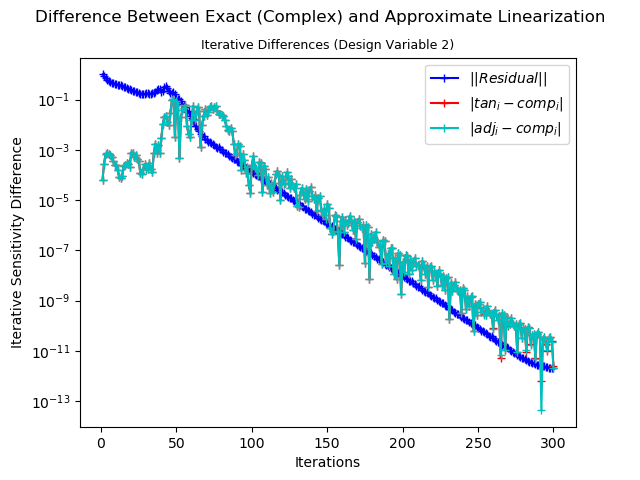}}
\end{subfigmatrix}
\caption{Difference between adjoint, tangent, and complex sensitivities at each iteration for a linear tolerance of $1e-1$}
\label{fig:2ndSensComp_LT1e-1}
\end{figure}

\begin{figure}
\begin{subfigmatrix}{2}
  \subfigure[Design Variable
1]{\includegraphics[keepaspectratio]{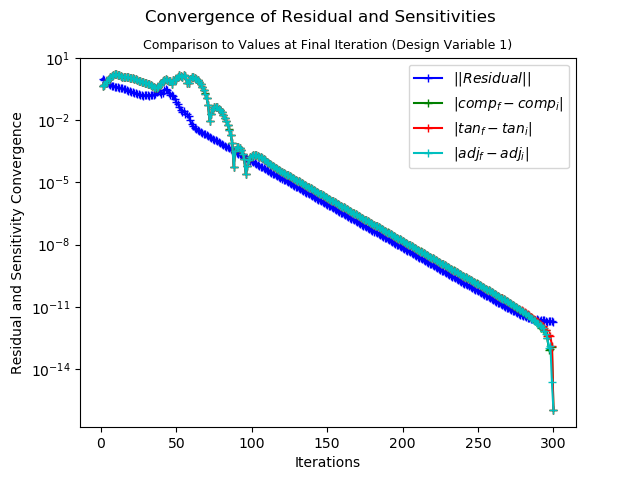}}
  \subfigure[Design Variable
2]{\includegraphics[keepaspectratio]{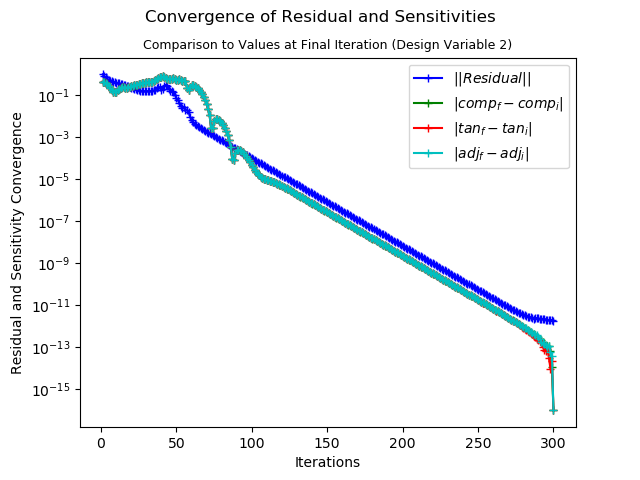}}
\end{subfigmatrix}
\caption{Convergence of nonlinear problem and sensitivities for linear tolerance of $1e-4$:
measured by $L_2$-norm of the residual and the difference between current and final sensitivity values}
\label{fig:2ndSensConv_LT1e-4}
\end{figure}

\begin{figure}
\begin{subfigmatrix}{2}
  \subfigure[Design Variable
1]{\includegraphics[keepaspectratio]{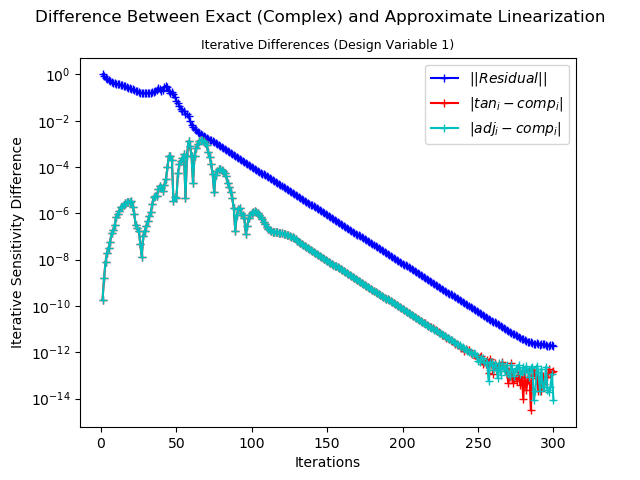}}
  \subfigure[Design Variable
2]{\includegraphics[keepaspectratio]{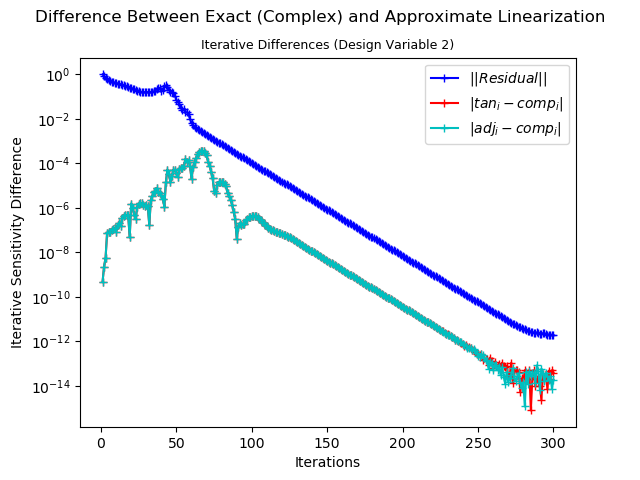}}
\end{subfigmatrix}
\caption{Difference between adjoint, tangent and complex sensitivities at each iteration for a linear tolerance of $1e-4$}
\label{fig:2ndSensComp_LT1e-4}
\end{figure}

This form is used as it is representative of a typical implicit solver used in
CFD problems, and the expectation is that we see similar behavior to that seen
in the previous section. In this case, rather than observing quadratic
convergence in both the primal and tangent/adjoint linearizations, we expect to
see linear convergence in the primal problem and its linearizations. To
demonstrate this we show
sister plots to those of the previous section, but for the inexact Newton runs.
Figures (\ref{fig:2ndSensConv_LT1e-1}, \ref{fig:2ndSensComp_LT1e-1}) show the
behavior for a linear tolerance of $1e-1$ and Figures
(\ref{fig:2ndSensConv_LT1e-4}, \ref{fig:2ndSensComp_LT1e-4}) show the behavior
for a tolerance of $1e-4$. These are the same tolerances we portrayed in the
previous section, and we can see the same expected convergence of the inexactly
linearized tangent and adjoint computed sensitivities to the complex linearization. We can see again, and as
hypothesized by the error bounds in this work, that the maximum iterative
difference is again on the order of the linear system tolerance, and that as the
analysis converges so do the adjoint, tangent, and complex sensitivities to each other,
down to nearly machine precision.

\begin{figure}
 \begin{subfigmatrix}{2}
    \centering
    \subfigure[Maximum Iterative Difference]{\includegraphics[keepaspectratio, width=.45\textwidth]{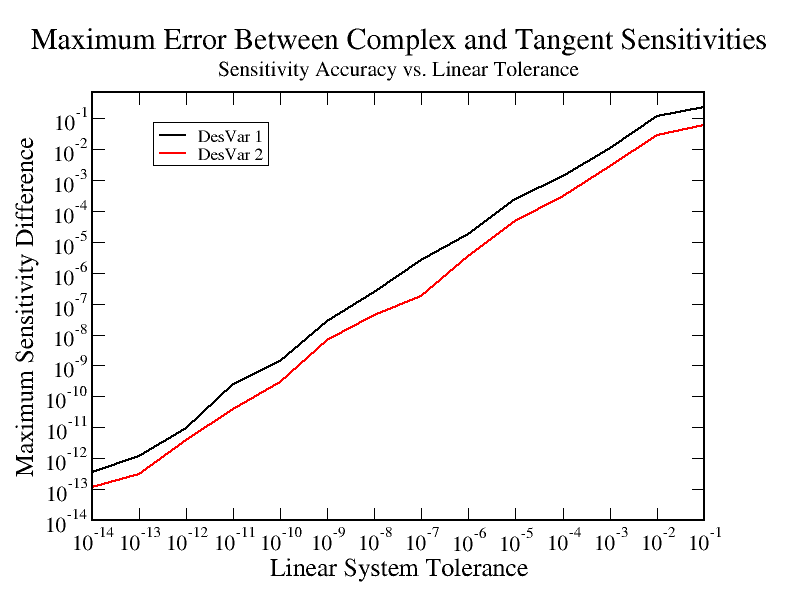}}
    \subfigure[Minimum Iterative Difference]{\includegraphics[keepaspectratio, width=.45\textwidth]{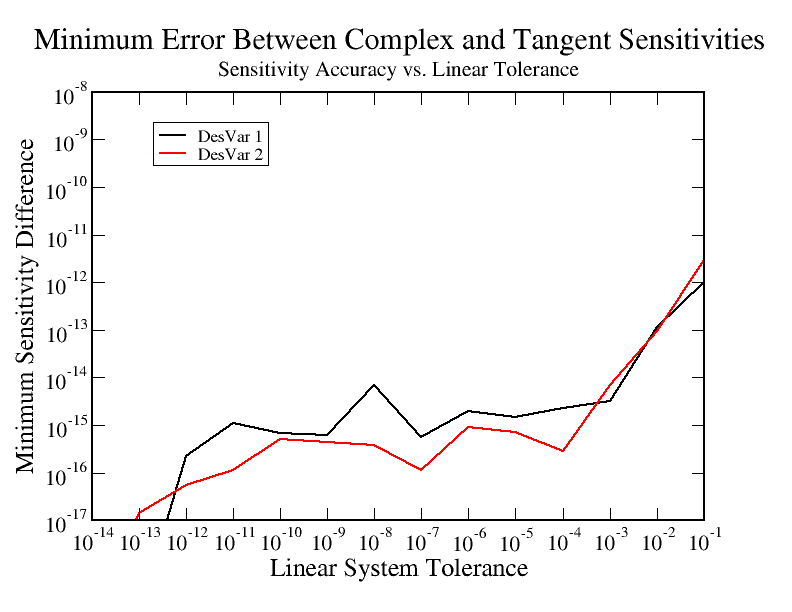}}
\end{subfigmatrix}   
\caption{Maximum and minimum difference between tangent and complex sensitivities as a function of linear tolerance}.
\label{fig:ErrConv_IQN}
\end{figure}

Having seen the impact of the tighter linear system tolerance on the maximum
iterative difference between the adjoint, tangent, and complex sensitivities in Figure (\ref{fig:ErrConv_QN}), and seeing
the expected behavior as shown in the theoretical error bound in Proposition (\ref{prop:TanErrIneq}) and the previous
section, we then simulate to check the impact of linear tolerances from $1e-1$
to $1e-14$. When we plot the maximum iterative difference in Figure
(\ref{fig:ErrConv_IQN}) we see the expected behavior over that parameter sweep.
This confirms the theoretical bound applies to a more general solver, one in which the
left hand side is not an exact linearization of the right hand side. Again, we see similar behavior between the tangent and adjoint linearizations but not identical as we have lost strict duality as we have not maintained the same number of linear iterations in the adjoint, tangent, and primal problems.

\subsection{Results for Inexactly Linearized Explicit Runge-Kutta Solver}

To check the effect of inexact linearization of the fixed point iteration for an explicit Runge-Kutta solver, we
calculate the flow gradients for the base stage, and then hold them frozen through
the fixed point iteration that solves the primal problem. For the linearization
of the fixed point iteration we still update the gradients at each step and
compute the Jacobian of the flow accordingly. To clarify the implementation we refer to the spatial residual as a summation over the edges or faces that construct the control volume, where the summation depends on the reconstructed states ($u_L, u_R$) given by

\begin{equation}
    R = \sum_{i=1}^{n_{edge}} F^{\bot}(u_L, u_R, n_{e_i}),
\end{equation}
Where $F^{\bot}$ and $n_{e_i}$ are defined as the inviscid flux normal to the face and the face unit normal.
The states $u_L$ and $u_R$ are reconstructed from the respective cell centers to the midpoint of the face separating cells j and k, and it becomes clearer how this inaccuracy is introduced. The reconstruction is implemented as
\begin{equation}
\begin{aligned}
    u_L &= u_j + \nabla u_j \cdot \vec{r}_{j}, \\
    u_R &= u_k + \nabla u_k \cdot \vec{r}_{k}, \\
\end{aligned}
\end{equation}
where $u_j, u_k$ are the cell center values at cells $j$ and $k$ respectively, and $\vec{r}_j, \vec{r}_k$ are the vectors from the respective element centers to the midpoint of the face separating the elements.
The explicit Runge-Kutta solver is expressed below, 
iterating $m$ through pseudo time
and $l$ through stages 1 to 5:

\begin{equation}
u^{m,l} = u^{m,0} + CFL\alpha^{l-1} \frac{\Delta t}{vol} R(u^{m,l-1}),
\end{equation}
with the end of the sub-stage time-stepping being governed as follows:
\begin{equation}
u^{m,0} = u^{m-1,5} .
\end{equation}

The primal problem fixed-point iteration uses the flow gradients from the base stage, i.e.:
\begin{equation}
\begin{aligned}
    u_L^l &= u_j^l + \nabla u_j^{m,0} \cdot \vec{r}_{j}, \\
    u_R^l &= u_k^l + \nabla u_k^{m,0} \cdot \vec{r}_{k}, \\
\end{aligned}
\end{equation}

but the inexact tangent and adjoint linearizations recompute the gradients at each stage as shown below.
\begin{equation}
\begin{aligned}
    u_L^l &= u_j^l + \nabla u_j^{m,l} \cdot \vec{r}_{j}, \\
    u_R^l &= u_k^l + \nabla u_k^{m,l} \cdot \vec{r}_{k}. \\
\end{aligned}
\end{equation}

\begin{figure}
 \begin{subfigmatrix}{2}
    \centering
    \subfigure[Design Variable 1]{\includegraphics[keepaspectratio, width=.45\textwidth]{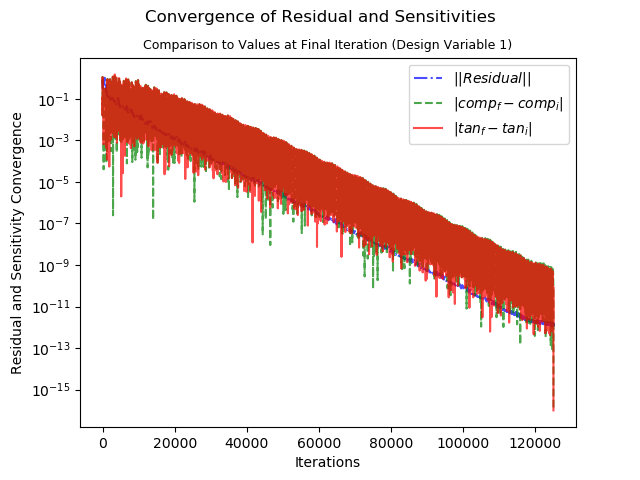}}
    \subfigure[Design Variable 2]{\includegraphics[keepaspectratio, width=.45\textwidth]{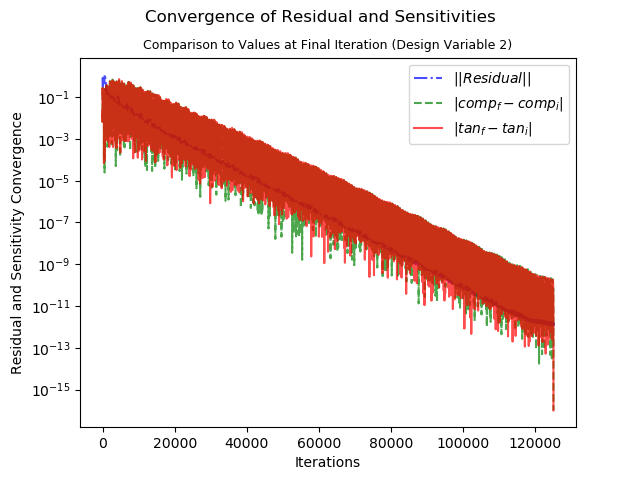}}
\end{subfigmatrix}   
\caption{Convergence of nonlinear problem and sensitivities for an inexactly linearized explicit Runge-Kutta solver:
measured by $L_2$-norm of the residual and the difference between current and final sensitivity values}
\label{fig:RKSensConv}
\end{figure}

This allows us to look at the Runge-Kutta scheme as some A operator depending on
the flow state and design variables that multiplies the 0-stage residual to get
$\delta u^m = u^{m,5} - u^{m,0} = A(u^{m,0})R(u^{m,0})$, where the linearization of A is inexact due to the
inconsistent handling of the gradients. When we look at the convergence of the
flow sensitivities to their respective final values as we did in the two
previous sections, Figure (\ref{fig:RKSensConv}) shows them converging to their
respective final values as the primal problem converges.

\begin{figure}
 \begin{subfigmatrix}{2}
    \centering
    \subfigure[Design Variable 1]{\includegraphics[keepaspectratio, width=.45\textwidth]{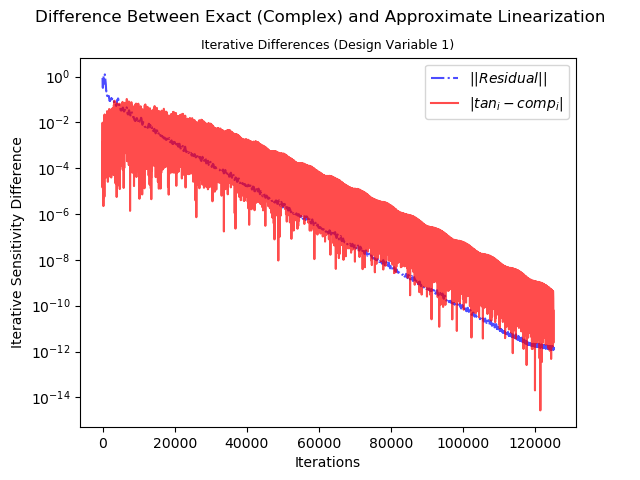}}
    \subfigure[Design Variable 2]{\includegraphics[keepaspectratio, width=.45\textwidth]{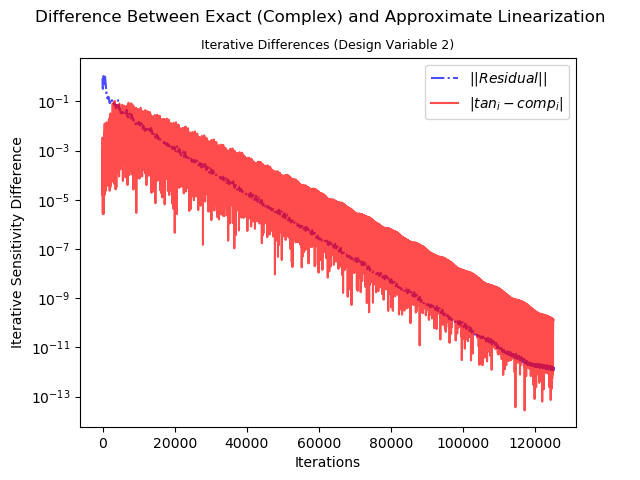}}
\end{subfigmatrix}   
\caption{Difference between tangent and complex sensitivities at each iteration for an inexactly linearized explicit Runge-Kutta solver}
\label{fig:RKSensComp}
\end{figure}

Looking at the difference between the complex (exact) linearization and the approximate tangent linearization portrayed in Figure (\ref{fig:RKSensComp}) shows the expected behavior. We see the error due to the inexact linearization of the Runge-Kutta scheme goes away at the rate of the primal problem convergence as shown in the proof section and the previous results. Due to the high number of iterations for the explicit Runge-Kutta solver a visualization of the comparison between the adjoint sensitivities and the complex-step finite difference computed ones is not possible, as this would be on the order of the number of iterations squared, but the expected behvaior is shown for the tangent as discussed previously, and the adjoint is the dual and has similar convergence properties.

\FloatBarrier
\section{Conclusion and Future Work}
\indent  The motivation for this work was to better understand
 the successful usage of the pseudo-time accurate adjoint algorithm when
 applied to sensitivity computation and optimization for nonconvergent
 simulations \cite{PTA-iqn}; to that end we sought to better understand the
 error introduced to the sensitivity computation by using inexact
 linearizations of fixed point iterations. In this paper, we showed that the
 error between the approximate linearization and the complex
 (exact) linearization of a fixed point iteration is a function of both the
 satisfaction of the discretized PDE denoted by the residual norm and the
 accuracy of the approximation of the fixed point iteration itself. We proved
 that in the limit of machine zero convergence of the nonlinear problem this
 approximation error vanishes. Furthermore, we have shown that we can obtain
 reasonable sensitivities for optimization without exact linearization of the
 fixed point iteration or full convergence of the non-linear problem, such that
 the error in the sensitivities depends on level of convergence in the
 non-linear and linear problems combined. This allows a user to
 select for the level of accuracy they desire in the sensitivities used for the
 optimization process. Finally, we demonstrated this behavior in a CFD code for
 an exact quasi-Newton solver, an inexact quasi-Newton solver and a five stage
 low-storage explicit Runge-Kutta solver. These proofs and results can also be
 useful in the analysis and implementation of automatic differentiation
 software, where the topic of error due to linearization of inexact linear
 solution processes is under investigation. Further future work in this area
 will relate to investigating other less accurate approximate linearizations
 and their use in optimization as well as using approximate linearizations in
 the aformentioned "piggy-back" iterations of the one-shot adjoint, where the
 authors believe similar proofs on error convergence can be shown and used to
 great effect. An additional application is that, having shown
 the convergence of the approximate linearizations parallels that of the
 non-linear problem, we can pair these approximate linearizations with work on
 design under inexact PDE constraints \cite{SKN} to get a more efficient design
 process for well-converging simulations.

\section{Acknowledgments}
This work was supported in part by NASA Grant NNX16AT23H and the NASA Graduate Aeronautics Scholars Program. Computing time was provided by ARCC on the Teton supercomputer.

\section{Conflict of Interest}
The authors declare that they have no conflict of interest.

\bibliography{references}
\bibliographystyle{spmpsci}
\end{document}